\documentclass[11pt,french]{smfart}

\usepackage[T1]{fontenc}
\usepackage[francais, english]{babel}

\usepackage{amssymb}

\usepackage{graphicx}
\usepackage{array}

\title[Sommation de Borel par s\'eries de factorielles]{Sommation effective
  d'une somme de Borel par s\'eries de factorielles}

\alttitle{Effective Borel-resummation by factorial series}
 
\author{Eric Delabaere}
\address{D\'epartement de Math\'ematiques, UMR CNRS 6093,
Universit\'e d'Angers, 2 Boulevard Lavoisier, 49045 Angers Cedex 01,
France.}
\email{eric.delabaere@univ-angers.fr}
\urladdr{}

\author{Jean-Marc Rasoamanana}
\address{D\'epartement de Math\'ematiques, UMR CNRS 6093,
Universit\'e d'Angers, 2 Boulevard Lavoisier, 49045 Angers Cedex 01,
France.}
\email{jean-marc.rasoamanana@univ-angers.fr}
\urladdr{}
 
\newtheorem{Theorem}{Th\'eor\`eme}[section]
\newtheorem{Proposition}{Proposition}[section]
\newtheorem{Lemma}{Lemme}[section]

\theoremstyle{definition}
 
\newtheorem{Definition}{D\'efinition}[section]
\newtheorem{Notation}{Notation}[section]

\newtheorem{rem}{Remarque}[section]

\begin{document}
 
\begin{abstract}
Nous abordons dans cet article la question de la sommation effective
d'une somme de Borel d'une s\'erie par la 
 s\'erie de factorielles associ\'ee. Notre approche fournit un
 contr\^ole de l'erreur entre la somme de Borel recherch\'ee et les sommes
 partielles de la s\'erie de factorielles. Nous g\'en\'eralisons ensuite
 cette m\'ethode au cadre des s\'eries de puissances fractionnaires,
 apr\`es avoir d\'emontr\'e un analogue d'un th\'eor\`eme de
 Nevanlinna de sommation de Borel fine pour ce cadre.
\end{abstract}

\begin{altabstract}
In this article, we consider the effective resummation of a Borel 
sum by its associated factorial series expansion. Our approach
provides concrete estimates for the remainder term when truncating 
this factorial series. We then generalize a theorem of Nevanlinna
which gives us the natural framework to extend the factorial series 
method for Borel-resummable fractional power series expansions.
\end{altabstract}

\altkeywords{Borel-resummation, factorial series.}
 
\maketitle

\section{Introduction}\label{sec8}

Le probl\`eme du calcul effectif d'une somme de Borel a d\'ej\`a une
longue histoire. Celui-ci peut \^etre fait par ce que
Poincar\'e appelait la ``m\'ethode des astronomes'' \cite{Poin87}, et
qui n'est autre que la m\'ethode de sommation au plus petit terme
que Stokes employait d\'ej\`a dans
son article fondateur de 1857 \cite{Sto57} et qui trouve naturellement
sa place dans le cadre Gevrey  
\cite{RamSch96, Can99}.

Bien d'autres m\'ethodes de sommation ont \'et\'e d\'evelopp\'ees
depuis. Ainsi, l'utilisation des
approximants de Pad\'e fait sont apparition d\`es les ann\'ees 1970 en
physique math\'ematique (voir, e.g., \cite{Sim82}), \`a la suite
notamment de la red\'ecouverte de la sommation de Borel par le 
``groupe de Saclay'' de physique
th\'eorique, 
avant d'\^etre d\'evelopp\'ee  
d'un point de vue algorithmique dans le cadre Gevrey \cite{Tho90}.  Un
point de vue diff\'erent, bas\'e sur l'utilisation de  transformations conformes,
fournit la m\'ethode expos\'ee dans \cite{Bal}.
Dans la mouvance des  id\'ees de Dingle \cite{Din73}, de Ecalle
 \cite{Ec81-1, Ec81-2, Ec85, CNP2, DP99, D94}, et sous
l'impulsion de Berry-Howls \cite{BeH91}, l'\'ecole anglo-saxonne
d'asymptotique exponentielle a quant \`a elle developp\'e  l'outil
de l'hyperasymptotique (voir, e.g.,  \cite{Old96, Old97, Old98}) pour lequel
la structure r\'esurgente  des
objets \`a sommer joue un r\^ole central \cite{Delab02, D06}. 
D'autres m\'ethodes effectives de sommation et  des applications en
  physique sont d\'etaill\'ees dans \cite{Jen}.
 
Nous allons nous pencher ici sur la m\'ethode des s\'eries de
factorielles. Cette m\'ethode, en th\'eorie exacte et non approch\'ee,
n'est pas nouvelle puisqu'elle 
remonte  \`a Watson \cite{Wn12}, Nevanlinna \cite{Nev18} et
N\"orlund \cite{Nor26}, voir aussi  \cite{Mal95, W65}. Notre apport
par rapport \`a la litt\'erature existante sur ce sujet est
double. D'une part, ayant en vue des m\'ethodes effectives, notre
objectif sera de fournir des estimations  de l'erreur
commise par sommation partielle  des s\'eries de factorielles. D'autre
part, nous proposerons une g\'en\'eralisation de la sommation par s\'eries
de factorielles au cadre des s\'eries de puissances fractionnaires
Borel sommables.

La structure de l'article est la suivante. 
La section \ref{somfact1} est consacr\'ee \`a un certain nombre de
rappels sur la sommation d'une s\'erie sommable de Borel par s\'eries
de factorielles. La section \ref{BorEffectif} fournit une approche
originale de la sommation par s\'eries de factorielles, d\'ebouchant
sur un contr\^ole effectif de l'erreur commise par sommation
partielle comme explicit\'e au  th\'eor\`eme \ref{thmsomfact1ter}. La
sommation de Borel de s\'eries de puissances fractionnaires est
abord\'ee en section \ref{somfact2}, l'objectif principal \'etant le
th\'eor\`eme \ref{thmNevfinram} de sommation de Borel fine. Celui-ci
fournit le cadre naturel pour la sommation par s\'eries de
factorielles g\'en\'eralis\'ee d\'evelopp\'ee en section
\ref{FactGene}, et son r\'esultat principal, le th\'eor\`eme
\ref{thmsomfact2bis}. Diff\'erents exemples illustrent les
proc\'edures de sommation effective.

Notons pour conclure cette introduction que cet article n'a pas pour
motivation de promouvoir telle m\'ethode effective de sommation 
au d\'epend de telle autre. Notre approche est ici de d\'efinir
les conditions d'applications de la m\'ethode par s\'eries de
factorielles, et donc ses limitations. Enfin, s'il nous semble que dans un
cadre r\'esurgent, et en
parall\`ele  avec l'hyperasymptotique, les m\'ethodes expos\'ees ici
devraient pouvoir d\'eboucher sur le calcul effectif des coefficients
de Stokes, cette question ne sera n\'eanmoins pas abord\'ee dans cet
article, faute pour les auteurs d'y avoir r\'efl\'echi suffisamment.

\section{R\'esultats classiques}\label{somfact1}

\begin{figure}[thp]
\begin{center}
\includegraphics[width=2.2in]{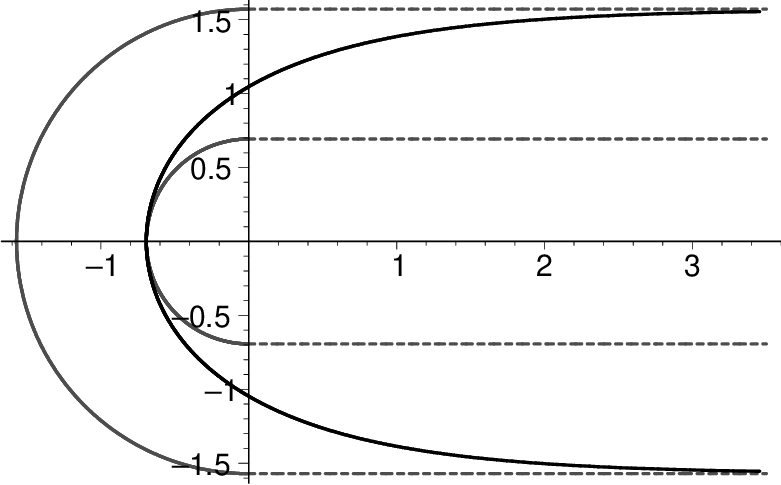} 
\caption{Les ouverts $\displaystyle  \mathcal{B}_{\ln(2)}
\subset \Delta  \subset \mathcal{B}_{\frac{\pi}{2}}$.
\label{fig:Som_Factorielle06}}
\end{center}
\end{figure}

\begin{Notation}
Dans tout l'article:
\begin{itemize}
\item  Pour $r>0$ et $\displaystyle \theta \in \frac{\mathbb{R}}{2 \pi \mathbb{Z}}$, 
 $\mathcal{B}_r (\theta)$ d\'esigne 
la bande ouverte 
$$\mathcal{B}_r (\theta) =\{\zeta \in
    \mathbb{C} \, /  \, d(\zeta,e^{i\theta} \mathbb{R}^+) < r\},$$
 o\`u  $d$ est la     distance euclidienne.\\
Pour $\theta =0$ on notera plus simplement $\mathcal{B}_r$ \`a la
place de  $\mathcal{B}_r (0)$.
\item On note $\Delta $ l'image du disque ouvert $D(1,1)$ 
de centre $1$ et de rayon $1$ par la transformation biholomorphe
$$s \in D(1,1) \mapsto \zeta = -\ln (s) \in \Delta.$$ 
 L'ouvert
$\Delta $ v\'erifie :
$$\displaystyle  \mathcal{B}_{\ln(2)}
\subset \Delta  \subset \mathcal{B}_{\frac{\pi}{2}} \hspace{5mm} 
\mbox{(cf. Fig. \ref{fig:Som_Factorielle06})}.$$
\item Pour tout $\lambda >0$, $\Delta_\lambda$ d\'esignera
  l'homoth\'etique de $\Delta$ d\'efini par:
$$\Delta_\lambda = \{\lambda \zeta \, /  \, \zeta \in  \Delta \}.$$
\end{itemize}
\end{Notation}

Nous r\'esumons maintenant la th\'eorie
classique de la sommation de Borel et de la sommation  par s\'eries de factorielles. Nous
renvoyons le lecteur \`a \cite{Mal95, W65} pour les
d\'emonstrations.

\subsection{Sommation de Borel}\label{rappelsomfact1}

Nous rappelons tout d'abord le th\'eor\`eme suivant de  Nevanlinna dit de 
``sommation de Borel fine'':

\begin{Theorem}\label{thmNevfin}
Soit $\displaystyle f(z)=\sum_{n=0}^{+\infty}\frac{a_n}{z^n} \in
\mathbb{C}[[z^{-1}]]_1$ une s\'erie Gevrey-1.
Notons $ \displaystyle \widetilde{f}(\zeta)=\sum_{n=1}^{+\infty}
\frac{a_n\zeta^{n-1}}{(n-1)!}  \in
\mathbb{C}\{\zeta\} $ sa transform\'ee de Borel.
Les assertions suivantes sont \'equivalentes~:
\begin{enumerate}
\item \label{item1} Il existe $r>0$ et $\displaystyle \theta \in \frac{\mathbb{R}}{2 \pi \mathbb{Z}}$ 
tels que $\widetilde{f}$ se prolonge
  analytiquement \`a l'ouvert $\mathcal{B}_r (\theta)$. 
De  plus,  il existe
    $A>0$, $B>0$ tels que pour tout $\zeta \in \mathcal{B}_r
    (\theta)$, $\displaystyle
    |\widetilde{f}(\zeta)| \leq Ae^{B|\zeta|}$. \\

\item Il existe $r>0$, $A >0$, $B>0$ et  une fonction  $\mbox{\sc  s}_{\theta} f(z)$
 holomorphe  dans $\{z \in \mathbb{C} / \, 
  \Re(ze^{i\theta})>B\}$  tels
  que,  pour $\Re(ze^{i\theta})> B $ et $n\geq 0$ : 
\begin{equation}\label{pluspetit}
\begin{array}{c}
 \displaystyle
  \Big|\mbox{\sc  s}_{\theta} f(z) -\sum_{k=0}^{n}
  \frac{a_k}{z^k}\Big| \leq R_{as}(r,A,B,n,ze^{i\theta})\\
\\
 \displaystyle R_{as}(r, A,B,n,z) = Ae^{Br}  \frac{ n!  }{r^n} \frac{1}{|z|^{n}(\Re
    (z)-B)}.
\end{array}
\end{equation}
\end{enumerate}
De plus, pour $\Re(ze^{i\theta})> B$, 
\begin{equation}\label{somborel}
\displaystyle
\mbox{\sc  s}_{\theta} f(z) =
a_0+\int_0^{\infty e^{i\theta}}\widetilde{f}(\zeta)e^{-z \zeta} \, d\zeta.
\end{equation}
\end{Theorem}

\begin{rem}\label{precpluspetit}
\begin{itemize}
\item L'hypoth\`ese  \ref{item1}) implique que pour tout $n \geq 0$,
\begin{equation}\label{plusfin}
|a_{n+1}| \leq  Ae^{Br} \frac{ n!  }{r^n} .
\end{equation}
\item Par l'in\'egalit\'e de Stirling il vient:
$$ R_{as}(r, A,B,n,z) < Ae^{Br}  \frac{ \sqrt{2\pi}
  n^{n+\frac{1}{2}}e^{-n+\frac{1}{12n}}}{r^n}  \frac{1}{|z|^{n}(\Re
    (z)-B)}.$$
Comme, \`a $z$ fix\'e, $\displaystyle \frac{n^n e^{-n}}{(r|z|)^n}$
atteint son minimum en $n=r|z|$, la {\em sommation au plus petit terme}
consiste \`a approcher $\mbox{\sc  s}_{\theta} f(z)$ par la somme 
$\displaystyle \sum_{k=0}^{n}
  \frac{a_k}{z^k}$ avec $n = \big[r|z|\big]$ o\`u $\big[ . \big]$ est la
partie enti\`ere.
%integer part en anglais
\end{itemize}
\end{rem}

\begin{Definition}
Dans le th\'eor\`eme (\ref{thmNevfin}), 
la fonction holomorphe $\mbox{\sc  s}_{\theta} f $ est la somme de Borel de $f$
relative \`a la direction $\theta$. 
\end{Definition}

\begin{rem}\label{somplus}
Le calcul de la somme de Borel $\mbox{\sc  s}_{\theta} f $ de $\displaystyle
f(z)=\sum_{n=0}^{+\infty}\frac{a_n}{z^n}$ relative \`a la direction
$\theta$ se ram\`ene au calcul de la somme de Borel 
$\mbox{\sc  s}_{0} f_\theta  $ de $\displaystyle
f_\theta (z)= f(ze^{-i\theta})  = \sum_{n=0}^{+\infty}\frac{a_n
  e^{in\theta} }{z^n}$ 
relative \`a la direction
$0$ par la relation: 
\begin{equation}\label{rotat}
\mbox{ pour tout } \Re (ze^{i\theta}) \gg 0, \hspace{3mm} 
\mbox{\sc  s}_{\theta} f (z) = \mbox{\sc  s}_{0} f_\theta (ze^{i\theta}).
\end{equation}
\end{rem}

\subsection{Sommation par s\'eries de factorielles}\label{factcla}

Nous consid\'erons maintenant la m\'ethode de calcul d'une
somme de Borel par s\'eries de factorielles. Suivant la remarque
\ref{somplus}, il suffit de se concentrer sur le calcul d'une somme de
Borel  relative \`a la direction
$\theta = 0$.

%%%%%%%%%%%%%%%%%%%%%%%%%%
Notre hypoth\`ese de travail est la suivante: la s\'erie Gevrey-1
$\displaystyle f(z)=\sum_{n=0}^{+\infty} \frac{a_n}{z^n}$ admet une 
 transform\'ee de Borel 
$\displaystyle
\widetilde{f}(\zeta)=\sum_{n=1}^{+\infty}\frac{a_n\zeta^{n-1}}{(n-1)!}$
qui  se prolonge holomorphiquement \`a l'ouvert
$\Delta$. On suppose par ailleurs:
\begin{equation}\label{continfty}
\displaystyle \exists A>0, \,\,  \exists B>0, \,\, \forall 
\zeta \in \Delta, \,\,  |\widetilde{f}(\zeta)|\leq Ae^{B |\zeta|}.
\end{equation}
%%%%%%%%%%%%%%%%%%%%%%%%
On notera que, puisque $\displaystyle  \mathcal{B}_{\ln(2)}
\subset \Delta$, l'hypoth\`ese (\ref{continfty}) induit que la somme
de Borel 
 $\mbox{\sc  s}_{0} f(z)$ de $f$ est  d\'efinie holomorphe pour
 $\Re(z) > B$.

L'application  $\zeta \mapsto s=e^{-\zeta}$ d\'efinissant  une
transformation biholomorphe  entre  l'ouvert $\Delta$ et le disque ouvert $D(1,1)$, nous pouvons
introduire: 
\begin{equation}\label{pourphi}
 \displaystyle
\phi(s)=\widetilde{f}(\zeta).
\end{equation}
L'application  $\phi$ est  holomorphe dans $D(1,1)$ et s'identifie
dans ce disque \`a la somme de sa s\'erie de Taylor:
\begin{equation}\label{coetay}
\displaystyle \phi(s)=\sum_{n=0}^{+\infty} b_n(1-s)^n, \hspace{5mm}
b_n = \frac{(-1)^n}{n !}  \frac{d^n \phi}{ds^n}(1).
\end{equation}

D'un point de vue formel, nous pouvons \'ecrire~:
$$ \begin{array}{lll}
\mbox{\sc  s}_{0} f(z) & = & \displaystyle a_0+\int_0^{+\infty}
\widetilde{f}(\zeta)e^{-z \zeta} \, d\zeta = 
a_0+\int_0^{+\infty} \phi(e^{-\zeta}) e^{-z \zeta} \, d\zeta 
 \\
& = & 
\displaystyle a_0+\sum_{n=0}^{+\infty}b_n\int_0^{+\infty}(1-e^{-\zeta})^ne^{-z\zeta}
\, d\zeta = a_0+\sum_{n=0}^{+\infty}b_n\int_0^1
(1-s)^ns^{z-1} \,ds \\
& = &\displaystyle a_0+ \sum_{n=0}^{+\infty} \frac{n! b_n}{z(z+1)\ldots(z+n)}.
\end{array}$$
C'est ce d\'eveloppement qui correspond \`a la sommation par  s\'erie de
  factorielles. Nous r\'ef\'erant \`a Malgrange \cite{Mal95}, 
la justification du calcul formel pr\'ec\'edent repose essentiellement sur le point
clef suivant:

\begin{Lemma}\label{lemmecruc}
Supposons qu'il existe $A>0$ et
$B>0$ tels que $\displaystyle \forall 
\zeta \in \Delta, \,\,  |\widetilde{f}(\zeta)|\leq Ae^{B |\zeta|}$. 
Alors pour tout $C > \max (B, 1)$ la s\'erie
$\displaystyle \sum_{n=1}^{+\infty}  \frac{|b_n|}{ n^C} $
converge.
\end{Lemma}

Ce lemme conduit alors au th\'eor\`eme suivant~:

\begin{Theorem}\label{thmsomfact1}
Avec les notations pr\'ec\'edentes, supposons qu'il existe $A>0$ et
$B>0$ tels que 
 $\displaystyle \forall 
\zeta \in \Delta, \,\,  |\widetilde{f}(\zeta)|\leq Ae^{B |\zeta|}$. \\
Alors la s\'erie de factorielles $\displaystyle a_0+\sum_{n=0}^{+\infty}
\frac{\Gamma(z) \Gamma(n+1) b_n}{\Gamma(z+n+1)}$ converge absolument
pour $\Re(z)> \max (B, 1)$
et repr\'esente la somme de Borel $\mbox{\sc  s}_{0} f(z)$ dans cet ouvert.
\end{Theorem}

L'utilisation de ce th\'eor\`eme n\'ecessite le calcul des
coefficients $b_n$  en fonction
des coefficients $a_n$ de la s\'erie formelle $f$. 
L'algorithme de Stirling \cite{Mal95} r\'epond \`a la question : 

\begin{Proposition}[Algorithme de Stirling]\label{algoStirling}
$$ \displaystyle \forall n \geq 0, \quad b_n=\frac{1}{n
  !}\sum_{k=1}^{n+1} (-1)^{n-k+1} \mathfrak{s}(n,k-1) a_k,$$ o\`u les
$\mathfrak{s}(n,k)$ sont les nombres de Stirling de premi\`ere esp\`ece.
\end{Proposition}

\begin{rem}
Notre d\'efinition des  nombres de Stirling de premi\`ere esp\`ece 
$\mathfrak{s}(n,k)$ est celle de \cite{Comtet}:  
%(Stirling cycle numbers, Stirling numbers of the first kind)  
$\displaystyle \prod_{k=0}^{n-1} (x-k) = \sum_{k=0}^n \mathfrak{s} (n,k) x^k$.
\end{rem}

Signalons pour m\'emoire que le th\'eor\`eme \ref{thmsomfact1} admet
une r\'eciproque, cf. \cite{Mal95}.

\section{Sommation de Borel effective}\label{BorEffectif}

Nous proposons \`a pr\'esent une
justification de la sommation par s\'erie de
factorielles. Diff\'erente des preuves classiques \cite{Mal95,
  W65}, son avantage r\'eside dans sa simplicit\'e 
et le fait qu'elle d\'ebouche sur un contr\^ole effectif. 

Notre hypoth\`ese de travail est celle du \S \ref{factcla}:  
la s\'erie Gevrey-1
$\displaystyle f(z)=\sum_{n=0}^{+\infty} \frac{a_n}{z^n}$ admet une 
 transform\'ee de Borel 
$\displaystyle \widetilde{f}(\zeta)$
qui  se prolonge holomorphiquement \`a l'ouvert
$\Delta$ et v\'erifie la condition (\ref{continfty}). Nous
consid\'erons de nouveau l'application $\phi$ d\'efinie par
(\ref{pourphi}) et sa s\'erie de Taylor (\ref{coetay}).

Tirons tout d'abord quelques cons\'equences de (\ref{continfty}) pour
les d\'eriv\'ees $\displaystyle \frac{d^{n}\phi}{ds^n}(s)$, $s \in ]0,1]$.
Par l'\'egalit\'e de Cauchy  et  pour $s \in ]0,1]$,
\begin{equation}\label{lab0}
\frac{d^{n}\phi}{ds^{n}}(s) = \frac{n!}{2i\pi}\oint
\frac{\phi(t)}{(t-s)^{n+1}}  dt = \frac{n!}{2i\pi s^{n+1}}\oint
\frac{\phi(t)}{(\frac{t}{s}-1)^{n+1}}  dt . 
\end{equation}
Dans (\ref{lab0}), l'int\'egration se fait le long d'un lacet dont un
param\'etrage possible est:  
\begin{equation}\label{lab1}
 t(\alpha) = s(1+re^{i\alpha}), \hspace{3mm} \alpha \in [0,2\pi],
\hspace{3mm} r \in ]0,1[ \mbox{ fix\'e}.
\end{equation}
Nous faisons \`a pr\'esent le changement de variable   $v \in D(1,1)  \mapsto
t=e^{-v} \in \Delta $. A $s \in ]0,1]$
correspond  $\zeta = -\ln(s) \in \mathbb{R}^+$ tandis qu'au lacet
(\ref{lab1})  est associ\'e le lacet
\begin{equation}\label{lab2}
 v(\alpha) = \zeta - \ln(1+re^{i\alpha}), \hspace{3mm} \alpha \in [0,2\pi],
\hspace{3mm}  0 < r <1.
\end{equation}
L'\'egalit\'e (\ref{lab0}) devient: pour tout $\zeta \in \mathbb{R}^+$
et tout $r \in ]0,1[$,
$$
\begin{array}{ll}
\displaystyle \frac{d^{n}\phi}{ds^{n}}(e^{-\zeta}) &  \displaystyle 
= -\frac{n! \, e^{\zeta(n+1)}}{2i\pi}\oint
\frac{\phi(e^{-v})}{(e^{-v+\zeta}-1)^{n+1}}e^{-v} dv \\
                                &  \displaystyle
= -\frac{n! \,   e^{\zeta(n+1)} }{2i\pi}\oint
\frac{\widetilde{f}(v)}{(e^{-v+\zeta}-1)^{n+1}}e^{-v} dv \\
                              & \displaystyle 
= \frac{n! \, e^{n \zeta} }{2\pi}\int_0^{2\pi}
\frac{\widetilde{f}( \zeta - \ln(1+re^{i\alpha})  )}{ (re^{i\alpha} )^n}
   d\alpha .
\end{array}
$$
Puisque pour $|\tau|<1$, $|\ln(1+\tau)| \leq - \ln(1-|\tau|)$
on en d\'eduit par (\ref{continfty}) que pour tout $\zeta \in \mathbb{R}^+$ et tout $r \in ]0,1[$,
\begin{equation}\label{merveille}
\begin{array}{ll}
\displaystyle \left| \frac{d^{n}\phi}{ds^{n}}(e^{-\zeta}) \right| & 
\displaystyle \leq  \frac{n!   \, e^{n \zeta}}{2\pi \, r^n }\int_0^{2\pi} 
Ae^{B(\zeta - \ln(1-r))}  \, d\alpha \\
              &  \displaystyle \leq   \frac{A \,n!  }{ r^n (1-r)^B}  \, e^{(n+B) \zeta}.
\end{array}
\end{equation}
Comme $r^{n} (1-r)^B$ atteint son
maximum en $\displaystyle r = \frac{n}{n+B}$ sur $]0,1[$, nous retiendrons
que:
\begin{equation}\label{merveillebis}
\forall \zeta \in \mathbb{R}^+, \, \, \,  \left|
  \frac{d^{n}\phi}{ds^{n}}(e^{-\zeta}) \right| \leq 
  \frac{A (n+B)^{n+B} \,n!  }{ B^B n^n}  \, e^{(n+B) \zeta}.
\end{equation}
Ceci \'etabli, consid\'erons maintenant la somme de Borel $\mbox{\sc
  s}_{0} f(z)$ de $f$.
Pour  tout $N \geq 0$ et $\Re (z) > B $, nous avons
\begin{equation}\label{diff1}
\mbox{\sc  s}_{0} f(z) - \Big( a_0+\sum_{n=0}^{N}
\frac{\Gamma(z) \Gamma(n+1) b_n}{\Gamma(z+n+1)} \Big) =
\int_0^{+\infty} \Big(\phi(e^{-\zeta}) - \sum_{n=0}^{N} b_n
(1-e^{-\zeta})^n \Big)e^{-z \zeta}   \,d\zeta .
\end{equation}
Maintenant par int\'egrations par parties, en utilisant (\ref{coetay})
et la majoration (\ref{merveillebis}) pour les questions d'int\'egrabilit\'e,
\begin{equation}
\begin{array}{c}
\displaystyle \mbox{\sc  s}_{0} f(z) - \Big( a_0+\sum_{n=0}^{N}
\frac{\Gamma(z) \Gamma(n+1) b_n}{\Gamma(z+n+1)} \Big) = \hspace{90mm} \,\\
\displaystyle \hspace{10mm}  \frac{(-1)^{N+1}}{z (z+1)\cdots (z+N)}\int_0^{+\infty}
\frac{d^{N+1}\phi}{ds^{N+1}}(e^{-\zeta}) e^{-(z+N+1) \zeta}
\,d\zeta .
\end{array}
\end{equation}
Par suite, en vertu de (\ref{merveillebis}): pour  tout $N \geq 0$ et $\Re (z) > B $,
\begin{equation}
\begin{array}{l}
\displaystyle \left| \mbox{\sc  s}_{0} f(z) - \Big( a_0+\sum_{n=0}^{N}
\frac{\Gamma(z) \Gamma(n+1) b_n}{\Gamma(z+n+1)} \Big)\right|\\ 
\\
\hspace{5mm}  \displaystyle \leq  \frac{A (N+B+1)^{N+B+1}  }{ B^B (N+1)^{N+1}} 
\frac{(N+1)! }{\left| z (z+1)\cdots (z+N) \right|}\int_0^{+\infty} e^{(B-
  \Re(z)) \zeta}    \,d\zeta\\
\\
\hspace{5mm} \displaystyle  \leq   \frac{A}{B^B}\frac{  \, (N+B+1)^{N+B+1}
  }{ (N+1)^{N} } 
\left| \frac{\Gamma (z) N! }{\Gamma(z+N+1) (\Re(z) -B)} \right|.
\end{array}
\end{equation}
Nous r\'esumons le r\'esultat obtenu:
\begin{Proposition}\label{thmsomfact1bis}
On  suppose que  la s\'erie 
$\displaystyle f(z)=\sum_{n=0}^{+\infty} \frac{a_n}{z^n} \in
\mathbb{C}[[z^{-1}]]_1$ admet une 
 transform\'ee de Borel 
$\displaystyle \widetilde{f}(\zeta)$
qui  se prolonge holomorphiquement \`a l'ouvert
$\Delta$, et qu'il existe   $A>0$ et
$B>0$ tels que pour tout $\zeta \in \Delta$,
$\displaystyle |\widetilde{f}(\zeta)|\leq Ae^{B |\zeta|}$.
Alors, pour  tout $N \geq 0$ et $\Re (z) > B $, 
\begin{equation}\label{remarquable}
\begin{array}{c}
\displaystyle \left| \mbox{\sc  s}_{0} f(z) - \Big( a_0+\sum_{n=0}^{N}
\frac{\Gamma(z) \Gamma(n+1) b_n}{\Gamma(z+n+1)} \Big)\right| \leq R_1(A,B,N,z)\\
\\
\displaystyle  R_1(A,B,N,z) =   \frac{A}{B^B}\frac{ (N+B+1)^{N+B+1}
  }{ (N+1)^{N} } 
\left| \frac{\Gamma (z) \Gamma (N+1)}{\Gamma(z+N+1) (\Re(z) -B)} \right|,
\end{array}
\end{equation}
o\`u $\mbox{\sc  s}_{0} f$ d\'esigne la somme de Borel de $f$, les
$b_n$ \'etant d\'eduits des $a_n$  par la proposition \ref{algoStirling}.
\end{Proposition}

La proposition  \ref{thmsomfact1bis} et le th\'eor\`eme \ref{thmsomfact1}  induisent
le r\'esultat suivant:

\begin{Theorem}\label{thmsomfact1ter}
On  suppose que  la  transform\'ee de Borel 
$\displaystyle \widetilde{f}(\zeta)$ de 
la s\'erie 
$\displaystyle f(z)=\sum_{n=0}^{+\infty} \frac{a_n}{z^n} \in
\mathbb{C}[[z^{-1}]]_1$ se prolonge holomorphiquement \`a l'ouvert
$\Delta_\lambda$, $\lambda>0$, et qu'il existe   $A>0$ et
$B>0$ tels que pour tout $\zeta \in \Delta_\lambda$,
$\displaystyle |\widetilde{f}(\zeta)|\leq Ae^{B|\zeta|}$. Alors: \\
$\bullet$   la s\'erie de factorielles $\displaystyle 
a_0  +\lambda \sum_{n=0}^{N}
\frac{\Gamma(\lambda z) \Gamma(n+1) b_n^{(\lambda)}}{\Gamma(\lambda
  z+n+1)}$ converge  absolument pour $\Re (z) > \max( B, 1/\lambda)$, de
somme $\mbox{\sc  s}_{0} f (z)$, o\`u $\mbox{\sc  s}_{0} f$ 
d\'esigne la somme de Borel de $f$, les
$b_n^{(\lambda)}$ \'etant d\'eduits des $a_n^{(\lambda)} = \lambda^{n-1} a_n$ 
 par la proposition \ref{algoStirling}.\\
$\bullet$ pour  tout $N \geq 0$ et $\Re (z) > B $, 
\begin{equation}\label{remarquablelamb}
\begin{array}{c}
\displaystyle 
 \left| \mbox{\sc  s}_{0} f (z) - \Big( a_0  +\lambda \sum_{n=0}^{N}
\frac{\Gamma(\lambda z) \Gamma(n+1) b_n^{(\lambda)}}{\Gamma(\lambda z+n+1)}
\Big)\right| \leq  R_{fact}(\lambda, A,B,N,z) \\
\\
\displaystyle R_{fact}(\lambda, A,B,N,z) = \hspace{80mm} \, \\
\displaystyle  \,\hspace{20mm}  \frac{A}{(\lambda B)^{\lambda B}} \frac{\, (N+\lambda
   B+1)^{N+\lambda B+1}
  }{ (N+1)^{N} } 
\left| \frac{\Gamma (\lambda z) \Gamma (N+1)}{\Gamma(\lambda z+N+1) (\Re(z) -B)}
\right|,
\end{array}
\end{equation}
\end{Theorem}

\begin{proof}
Posons  $\displaystyle \widetilde{f}_\lambda (\zeta) = \widetilde{f}(\lambda
\zeta)$, de sorte que
$\widetilde{f}_\lambda$ se prolonge  holomorphiquement  sur 
$\Delta$, et
 $ \displaystyle \forall \zeta \in \Delta$, 
$ \,|\widetilde{f}_\lambda(\zeta)| \leq Ae^{\lambda B|\zeta|}$. La fonction
 $\widetilde{f}_\lambda$ n'est autre que la transform\'ee de Borel de la
s\'erie formelle $\displaystyle f_\lambda(z) = \frac{1}{\lambda} f\left(
  \frac{z}{\lambda} \right) =  \sum_{n=0}^{+\infty}
\frac{a_n^{(\lambda)}}{z^n}$ avec $\displaystyle
a_n^{(\lambda)} = \lambda^{n-1} a_n$. 
Nous d\'eduisons du th\'eor\`eme  \ref{thmsomfact1}
que la s\'erie de factorielles $\displaystyle  a_0^{(\lambda)}+\sum_{n=0}^{N}
\frac{\Gamma(z) \Gamma(n+1) b_n^{(\lambda)}}{\Gamma(z+n+1)}$ converge
(absolument) vers $\mbox{\sc  s}_{0} f_\lambda (z)$ 
 pour $\Re (z) > \max( \lambda B, 1)$, et par la proposition 
 \ref{thmsomfact1bis} que
 pour  tout $N \geq 0$ et $\Re (z) > \lambda B $, 
$$
 \left| \mbox{\sc  s}_{0} f_\lambda (z) - \Big( a_0^{(\lambda)}+\sum_{n=0}^{N}
\frac{\Gamma(z) \Gamma(n+1) b_n^{(\lambda)}}{\Gamma(z+n+1)}
\Big)\right| \leq 
R_1 (A,\lambda B,N,z).
$$
o\`u les $b_n^{(\lambda)}$ sont d\'eduits des $a_n^{(\lambda)}$  
par l'algorithme de Stirling (proposition \ref{algoStirling}).\\
Par suite la s\'erie de factorielles $\displaystyle 
a_0  +\lambda \sum_{n=0}^{N}
\frac{\Gamma(\lambda z) \Gamma(n+1) b_n^{(\lambda)}}{\Gamma(\lambda
  z+n+1)}$ converge  (absolument) vers $\mbox{\sc  s}_{0} f
(z)$ pour $\Re (z) > \max( B, 1/\lambda)$, et 
pour  tout $N \geq 0$ et $\Re (z) > B $,
$$
 \left| \mbox{\sc  s}_{0} f (z) - \Big( a_0  +\lambda \sum_{n=0}^{N}
\frac{\Gamma(\lambda z) \Gamma(n+1) b_n^{(\lambda)}}{\Gamma(\lambda z+n+1)}
\Big)\right| \leq 
\lambda R_1 (A,\lambda B,N,\lambda z).
$$
\end{proof}

Le lemme suivant est une cons\'equence facile de la formule de Stirling.

\begin{Lemma}\label{estima}
Avec les notations du th\'eor\`eme \ref{thmsomfact1ter}, pour $\Re(z) >B$,
$$R_{fact}(\lambda, A,B,N,z)  \sim \frac{A e^{\lambda
    B\big(1-\ln(\lambda B)\big)}}{N^{\lambda \big( \Re(z)-B \big) -1}} \frac{ \left|  \Gamma
    (\lambda z)  \right| }{\Re(z) -B}$$
quand $N \rightarrow +\infty$.
\end{Lemma}

De la majoration
(\ref{merveillebis}) (pour $\lambda =1$) et de la relation
(\ref{coetay}) il  d\'ecoule:

\begin{Lemma}\label{estima2}
Avec les notations du th\'eor\`eme \ref{thmsomfact1ter}, pour tout $n
\geq 0$,
\begin{equation}\label{pourbn}
|b_n^{(\lambda)} | \leq  \frac{A (n+\lambda B)^{n+ \lambda B} }{
  (\lambda B)^{\lambda B} n^n}.
\end{equation}
\end{Lemma}

\begin{rem}
Le lemme \ref{estima}  d\'emontre la convergence de la s\'erie
de factorielles  $\displaystyle  a_0  +\lambda \sum_{n=0}^{+\infty}
\frac{\Gamma(\lambda z) \Gamma(n+1) b_n^{(\lambda)}}{\Gamma(\lambda
  z+n+1)}$ 
 pour $\displaystyle \Re(z)>B+ \frac{1}{\lambda}$, $B>0$,
celle-ci repr\'esentant  la somme de Borel $\mbox{\sc  s}_{0} f(z)$
dans cet ouvert.\\
Ce r\'esultat est plus faible que celui donn\'e par le th\'eor\`eme
\ref{thmsomfact1ter}. Le d\'ecalage s'explique par la majoration
obtenue au lemme \ref{estima2}, $\displaystyle |b_n| \leq \frac{A (n+B)^{n+B} }{ B^B
  n^n}$ (pour $\lambda =1$). 
Or  $\displaystyle \frac{A (n+B)^{n+B} }{ B^B
  n^n} \sim \frac{A e^B }{ B^B } n^B$,  \`a comparer avec le lemme
\ref{lemmecruc}.
\end{rem}
$\,$

\subsection{Un \'el\'ement de comparaison}

\begin{figure}[thp]
\begin{center}
\includegraphics[width=2.5in]{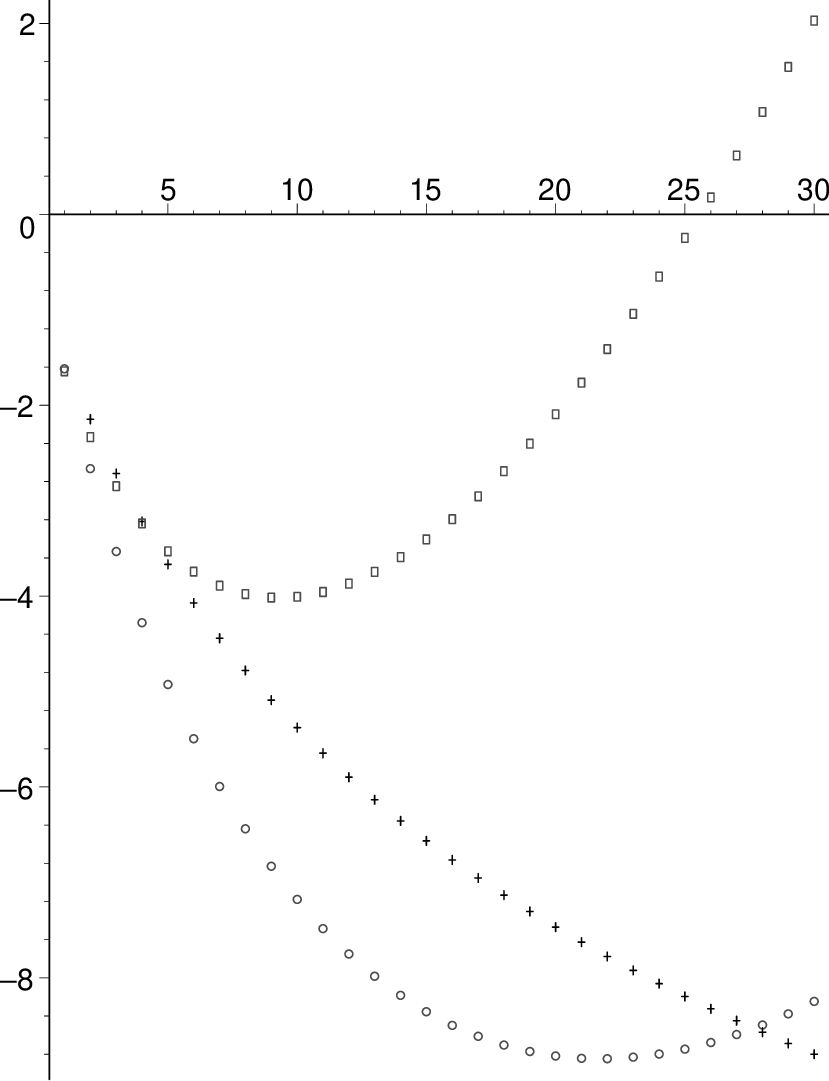} 
\caption{
Nous avons repr\'esent\'e pour $A=B=1$, $z=10+10i$ et 
 $n=0,\cdots, 30$, par $\Box$ les
  points de coordonn\'ees $[n, \log \left| R_{as}(\ln(2), 1,1,n,z)
  \right|]$,  par $\circ$ les
  points de coordonn\'ees $[n, \log \left| R_{as}(\frac{\pi}{2},
    1,1,n,z) 
  \right|]$, par $+$ les
  points de coordonn\'ees $[n, \log \left|R_{fact}(1 , 1, 1,n, z
  \right|]$.
\label{fig:Som_Factorielle07}}
\end{center}
\end{figure}

Au vu de la relation $\displaystyle  \mathcal{B}_{\ln(2)}
\subset \Delta  \subset \mathcal{B}_{\frac{\pi}{2}}$, nous allons  comparer
les estimations des restes $R_{as}(\ln(2), A,B,N,z)$, $\displaystyle
R_{as}(\frac{\pi}{2}, A,B,N,z)$ donn\'ees par (\ref{pluspetit}), 
et celle du reste $R_{fact}(1 , A,B,N+1,z)$ d\'ecrit par le
th\'eor\`eme \ref{thmsomfact1ter}. Nous prendrons $A=1$, $B=1$,
$z=10+10i$. Le r\'esultat est d\'ecrit par la figure \ref{fig:Som_Factorielle07}.

Des exemples d'applications seront donn\'es dans les sections qui suivent.

\section{Sommation de Borel de s\'eries de 
puissances fractionnaires}\label{somfact2}

La section \ref{somfact1} s'occupait de sommes de Borel de s\'eries
non ramifi\'ees. Nous allons \`a pr\'esent nous pencher sur le
probl\`eme de la sommation effective d'une somme de Borel  d'une
s\'erie de  puissances fractionnaires. Cette section a pour but de
d\'evelopper un analogue du th\'eor\`eme \ref{thmNevfin} de "sommation
de Borel fine".

\begin{Notation}
Par la suite $m \in \mathbb{N}^*$. Nous noterons par 
$\displaystyle 
\begin{array}{c}
\mathbb{C}_m \\
\pi \downarrow \\
\mathbb{C}
\end{array}$ \Big({\em resp.} $\displaystyle 
\begin{array}{c}
\mathbb{C}_m^\star \\
\pi \downarrow \\
\mathbb{C}^\star
\end{array}$\Big)
la surface de Riemann ramifi\'ee ({\em resp.} la surface de Riemann) 
\`a $m$ feuillets de $X^{\frac{1}{m}}$. 
\begin{itemize}
\item Nous identifierons l'\'el\'ement $x \in \mathbb{C}_m^\star$ au  couple 
$(|x|, \arg (x))$ o\`u $|x| \in \mathbb{R}^{+\star}$ et $\displaystyle
\arg (x) \in \frac{\mathbb{R}}{2 \pi m \mathbb{Z}}$ et on notera $x = |x|e^{i\arg(x)}$ :
$$\begin{array}{ccc}
x = r e^{i\arg(\theta)} \in \mathbb{C}_m^\star & \longleftrightarrow & 
(r, \theta) = \big(|x|, \arg (x)\big) \in
\mathbb{R}^{+\star} \times \frac{\mathbb{R}}{2 \pi m \mathbb{Z}}\\
\pi \downarrow  & & \downarrow \tilde \pi \\
\dot x = r e^{i\arg(\dot \theta)} \in \mathbb{C}^\star    &
\longleftrightarrow   &  
(r, \dot \theta) = \big( |\dot x|, \arg (\dot x) \big) \in
\mathbb{R}^{+\star} \times \frac{\mathbb{R}}{2 \pi \mathbb{Z}}
\end{array}
$$
\item  Pour $x \in \mathbb{C}_m^\star$ et
$k \in \mathbb{Z}$, on notera $x^{k/m}$ l'\'el\'ement de
$\mathbb{C}_m^\star$ de module $|x^{k/m}| = |x|^{k/m}$ et d'argument 
$\arg (x^{k/m}) = \frac{k}{m} \arg (x)$. 
\item  Pour $r>0$ et $\displaystyle \dot \theta \in  \frac{\mathbb{R}}{2 \pi \mathbb{Z}}$ on notera
$$\mathcal{D}_r (\dot \theta) = \pi ^{-1} \Big( \mathcal{B}_r (\dot
\theta) \Big) \subset \mathbb{C}_m, \hspace{5mm} 
\mathcal{D}_r^\star (\dot \theta) = \pi ^{-1} \Big( \mathcal{B}_r^\star  (\dot
\theta) \Big) \subset \mathbb{C}_m^\star .
$$
Pour $\dot \theta =0$ on \'ecrira plus simplement 
$\mathcal{D}_r  = \pi ^{-1} ( \mathcal{B}_r  )$ et $\mathcal{D}_r^\star  = \pi ^{-1} ( \mathcal{B}_r^\star  )$.
\item On note $\displaystyle \Omega = \pi ^{-1} ( \Delta )$  et pour 
tout $\lambda >0$, $\displaystyle \Omega_\lambda = \pi ^{-1} (
\Delta_\lambda )$. On d\'efinit $\Omega^\star$ et
$\Omega_\lambda^\star$ de fa\c con analogue.
\item Pour $B >0$, on notera $P(B)$ le demi-plan ouvert
$$P(B) = \{z \in  \mathbb{C}_m^\star \, / \,
|\arg(z)| \leq \frac{\pi}{2}, \, \Re
(\dot z ) > B\}.$$
\end{itemize}
\end{Notation}

\subsection{Somme de Borel}

Rappelons quelques d\'efinitions et propri\'et\'es \'el\'ementaires.

\begin{Definition}\label{Deframbo}
Soit $\displaystyle f(z)=\sum_{n=0}^{+\infty}
\frac{a_n}{z^{\frac{n}{m}}} \in
\mathbb{C}[[z^{-\frac{1}{m}}]]$ et $\displaystyle \widetilde{f}(\zeta)=\sum_{n=1}^{+\infty}
\frac{a_n \zeta^{\frac{n}{m}-1} }{\Gamma\Big(\frac{n}{m}\Big)} $ sa
transform\'ee de Borel. On suppose qu'il existe $r>0$ tel que
$\displaystyle \widetilde{f}(\zeta)$ d\'efinisse une fonction
holomorphe sur $\pi^{-1} \big( D(0,r)^\star \big)$,  $ D(0,r)^\star
= D(0,r) \backslash \{0\}$, o\`u 
$D(0,r)$ est le disque ouvert de centre $0$ de rayon $r$.\\
Soit $\theta \in \frac{\mathbb{R}}{2 \pi m
  \mathbb{Z}}$. On suppose qu'il existe un secteur ouvert $\Sigma (\theta,
\varepsilon) = \{x \in \mathbb{C}_m^\star, \, \, \arg (x) \in ]\theta
-\varepsilon, \theta + \varepsilon[ \}$, $\varepsilon >0$, tel que
$\displaystyle \widetilde{f}(\zeta)$ se prolonge analytiquement dans
ce secteur et:
\begin{equation}\label{majsect}
\exists \, A >0, \,\exists \, B >0, \, \forall \zeta \in \Sigma (\theta,
\varepsilon), \, \, |\widetilde{f}(\zeta) \zeta^\frac{m-1}{m} |\leq A e
^{B|\zeta|}.
\end{equation}
Alors la fonction 
\begin{equation}\label{somborelram}
\displaystyle
\mbox{\sc  s}_{\theta} f(z) =
a_0+\int_0^{\infty e^{i\theta}}\widetilde{f}(\zeta)e^{-z \zeta} \, d\zeta
\end{equation}
d\'efinie holomorphe dans le demi-plan ouvert $\displaystyle \{z \in  \mathbb{C}_m^\star \, / \,
|\arg(z)+\theta| \leq \frac{\pi}{2}, \, \Re
(\dot z e^{i \dot \theta}) > B\}$, est la somme de Borel de $f$ dans
la direction $\theta$.
\end{Definition}

\begin{rem}\label{somplusram}
Le calcul de la somme de Borel $\mbox{\sc  s}_{\theta} f $ de $\displaystyle
 f(z)=\sum_{n=0}^{+\infty}
\frac{a_n}{z^{\frac{n}{m}}}$ relative \`a la direction
$\theta \in \frac{\mathbb{R}}{2 \pi m
  \mathbb{Z}}$ se ram\`ene au calcul de la somme de Borel 
$\mbox{\sc  s}_{0} f_\theta  $ de $\displaystyle
f_\theta (z)= f(ze^{-i\theta})  = \sum_{n=0}^{+\infty}
\frac{a_n e^{-in\theta/m}}{z^{\frac{n}{m}}}$ 
relative \`a la direction
$0$ par la relation: 
$$
\mbox{\sc  s}_{\theta} f (z) = \mbox{\sc  s}_{0} f_\theta (ze^{i\theta}).
$$
\end{rem}

La remarque \ref{somplusram} permet de se limiter \`a l'\'etude de la
sommation de Borel pour la direction $\theta =0$, ce que nous ferons
dans la suite.

\begin{Proposition}\label{MagGevrey}
On se place dans le cadre de la d\'efinition \ref{Deframbo} avec
$\theta=0$. Soit $0 < \delta  < \pi/2$ et $\mu >1$. On note
$$P_{\delta, \mu} (B) = \{z \in  \mathbb{C}_m^\star \, / \,
|\arg(z)| \leq   \frac{\pi}{2} -\delta, \, |z| \geq 
\frac{\mu B}{\sin(\delta)}\}.$$
Alors il existe $C>0$ tel que 
\begin{equation}\label{majGevrey}
\forall \, z \in P_{\delta, \mu}  (B), \, \forall \, N \geq 0, \, 
\left| \mbox{\sc  s}_{0} f (z)  - \sum_{n=0}^{N}
\frac{a_n}{z^{\frac{n}{m}}}\right|  \leq  C^{1+N/m} \frac{\Gamma\left(
 1+ \frac{N}{m}\right)}{|z|^{1+N/m} }.
\end{equation}
\end{Proposition}

\begin{proof}
Quitte \`a raisonner avec $\displaystyle f(z) - \sum_{n=0}^{m-1}
\frac{a_n}{z^{\frac{n}{m}}}$ on peut supposer que $a_j=0$, $j=0,
\cdots, m-1$.
Soit  $0< b < r$ et $N \geq m$. Pour  $z \in P_{\delta, \mu} (B)$,
$$ \mbox{\sc  s}_{0} f (z) = 
\int_0^{\infty} \widetilde{f}(\zeta)e^{-z \zeta} \, d\zeta = 
\int_0^{b}\widetilde{f}(\zeta)e^{-z \zeta} \, d\zeta +
\int_b^{\infty}\widetilde{f}(\zeta)e^{-z \zeta} \, d\zeta.
$$
Par l'holomorphie de $\displaystyle \widetilde{f}(\zeta)$ sur
$\pi^{-1} \big( D(0,r)^\star \big)$ on peut \'ecrire
$$
\int_0^{b}\widetilde{f}(\zeta)e^{-z \zeta} \, d\zeta = \sum_{n=m}^{+\infty}
\int_0^{b} \frac{a_n \zeta^{\frac{n}{m}-1}
}{\Gamma\big(\frac{n}{m}\big)}e^{-z \zeta}\, d\zeta, $$
puis
$$\displaystyle 
\begin{array}{l}
\displaystyle \int_0^{b}\widetilde{f}(\zeta)e^{-z \zeta} \, d\zeta - \sum_{n=m}^{N}
\frac{a_n}{z^{\frac{n}{m}}} = \\
\displaystyle \hspace{25mm} - \sum_{n=m}^{N}
\int_b^{\infty} \frac{a_n \zeta^{\frac{n}{m}-1}
}{\Gamma\big(\frac{n}{m}\big)}e^{-z \zeta}\, d\zeta + 
\sum_{n=N+1}^{+\infty}
\int_0^{b} \frac{a_n \zeta^{\frac{n}{m}-1}
}{\Gamma\big(\frac{n}{m}\big)}e^{-z \zeta}\, d\zeta.
\end{array}
$$
En posant $\zeta = b t$ on obtient:
$$
\begin{array}{l}
\displaystyle \int_0^{b}\widetilde{f}(\zeta)e^{-z \zeta} \, d\zeta - \sum_{n=m}^{N}
\frac{a_n}{z^{\frac{n}{m}}} = \\
\displaystyle \hspace{25mm}  - \sum_{n=m}^{N} \frac{a_n b^{\frac{n}{m}}
}{\Gamma\big(\frac{n}{m}\big)}
\int_1^{\infty}  t^{\frac{n}{m}-1} e^{-z b t}\, dt
+ 
\sum_{n=N+1}^{\infty} \frac{a_n b^{\frac{n}{m}}
}{\Gamma\big(\frac{n}{m}\big)}
\int_0^{1}  t^{\frac{n}{m}-1} e^{-z b t}\, dt.
\end{array}
$$
Dans chacune des int\'egrales on peut majorer $\displaystyle
t^{\frac{n}{m}-1}$ par $\displaystyle
t^{\frac{N}{m}}$. On obtient ainsi:
$$\left| \int_0^{b}\widetilde{f}(\zeta)e^{-z \zeta} \, d\zeta - \sum_{n=m}^{N}
\frac{a_n}{z^{\frac{n}{m}}}\right|  \leq  
\left( \sum_{n=m}^{\infty } \frac{|a_n| b^{\frac{n}{m}}
}{\Gamma\big(\frac{n}{m}\big)} \right) \frac{\Gamma\left(
  \frac{N}{m}+1\right)}{b^{1+N/m} (\Re \dot z)^{1+N/m} }.
$$
Comme $z \in P_{\delta, \mu} (B)$ implique que 
$\Re(\dot z) \geq \sin(\delta)|z| \geq \mu B$, on en d\'eduit l'existence
d'une constante $c >0$ tel que pour tout $z \in P_{\delta, \mu} (B)$,
$$\left| \int_0^{b}\widetilde{f}(\zeta)e^{-z \zeta} \, d\zeta - \sum_{n=m}^{N}
\frac{a_n}{z^{\frac{n}{m}}}\right|   
\leq  c^{1+N/m} \frac{\Gamma\left(
 1+ \frac{N}{m}\right)}{|z|^{1+N/m} }.
$$
Par ailleurs l'hypoth\`ese (\ref{majsect}) implique :
$$\exists \, A >0, \,\exists \, B >0, \, \forall \zeta \in \Sigma (0,
\varepsilon), \, \, |\widetilde{f}(\zeta)  |\leq A e
^{B|\zeta|}.$$
Par suite, puisque  
$\Re(\dot z) - B \geq  (1-\frac{1}{\mu})\sin(\delta)|z| \geq (\mu-1)B$ 
pour tout $z \in P_{\delta, \mu} (B)$,
$$\left| \int_b^{\infty}\widetilde{f}(\zeta)e^{-z \zeta} \, d\zeta
\right| \leq \frac{e^{(B - \Re(\dot z))b} }{\Re (\dot z) -B} \leq 
\frac{e^{-(1-\frac{1}{\mu}) b \sin(\delta)|z|} }{(\mu-1)B}.
$$
Or, pour $\alpha >0$, $\displaystyle 
e^{-(1-\frac{1}{\mu}) b \sin(\delta)|z|}|z|^\alpha$ est maximal pour
$\displaystyle  |z| = \frac{\alpha}{(1-\frac{1}{\mu}) b
  \sin(\delta)}$. Donc,
$$\left| \int_b^{\infty}\widetilde{f}(\zeta)e^{-z \zeta} \, d\zeta
\right| \leq 
\frac{e^{-\alpha} }{(\mu-1)B} \left(\frac{\alpha}{(1-\frac{1}{\mu}) b
  \sin(\delta)}\right)^\alpha |z|^{-\alpha} .
$$
Le r\'esultat annonc\'e s'en d\'eduit alors par l'in\'egalit\'e de Stirling et par
l'in\'egalit\'e triangulaire.
\end{proof}

\begin{rem}\label{Conssect}
L'hypoth\`ese  (\ref{majsect}) faite dans la d\'efinition
\ref{Deframbo},  valable dans un secteur, induit que les
estimations Gevrey (\ref{majGevrey}) de la proposition \ref{MagGevrey}
s'\'etendent \`a un secteur d'ouverture plus grande que $\pi$. Ceci
implique l'unicit\'e de la somme de Borel.
\end{rem}

\subsection{Sommation de Borel fine} 

Nous allons \`a pr\'esent modifier quelque peu nos hypoth\`eses afin de
d\'emontrer un analogue du th\'eor\`eme \ref{thmNevfin}. Par la
remarque \ref{somplusram} il est loisible de se limiter \`a l'\'etude de la
sommation de Borel pour la direction $\theta =0$.

\begin{Theorem}\label{thmNevfinram}
Soit $\displaystyle f(z)=\sum_{n=0}^{+\infty}
\frac{a_n}{z^{\frac{n}{m}}} \in
\mathbb{C}[[z^{-\frac{1}{m}}]]_1$ une s\'erie Gevrey-1 
et $\displaystyle \widetilde{f}(\zeta)=\sum_{n=1}^{+\infty}
\frac{a_n \zeta^{\frac{n}{m}-1} }{\Gamma\Big(\frac{n}{m}\Big)}  \in
\zeta^{-1+\frac{1}{m}}\mathbb{C}\{\zeta^\frac{1}{m}\} $ sa
transform\'ee de Borel.
Les assertions suivantes sont \'equivalentes~:
\begin{enumerate}
\item  Il existe $r>0$ 
tel que $\widetilde{f}$ se prolonge
  analytiquement \`a l'ouvert $\mathcal{D}_r^\star$ de $\mathbb{C}_m$. 
De  plus,  
\begin{equation}\label{majDram}
\exists \, A >0, \,\exists \, B >0, \, \forall \zeta \in
\mathcal{D}_r^\star, \, \, 
|\widetilde{f}(\zeta) \zeta^\frac{m-1}{m} |\leq A e^{B|\zeta|}.
\end{equation}
\item Il existe $r>0$, $B>0$, $A_l >0$ pour $1\leq l\leq m$ et des fonctions
$\mbox{\sc  s}_{0} f_l(z)$
 holomorphe  dans $\Re(\dot z)>B$  tels
  que,  pour $\Re(\dot z)>B$, $n\geq 1$ et  $1\leq l\leq m$:
\begin{equation}\label{pluspetitram}
\begin{array}{c}
 \displaystyle
   \Big|\mbox{\sc  s}_{0} f_l(\dot z) -\sum_{j=1}^{n}
  \frac{a_{l,j}}{{\dot z}^{j}}\Big| \leq R_{as}(r,A_l,B,n, \dot z)\\
\\
 \displaystyle  \mbox{o\`u} \hspace{5mm}  f_l(\dot z)=\sum_{j=1}^{+\infty}
\frac{a_{l,j}}{{\dot z}^{j}}, \, \, a_{l,j} = a_{l+m(j-1)}.
\end{array}
\end{equation}
\end{enumerate}
De plus, pour $\Re (\dot z) >B$, 
$\displaystyle
\mbox{\sc  s}_{0} f_l(\dot z)  =
\int_0^{\infty}\widetilde{f_l}(\zeta)e^{-\dot z \zeta} \, d\zeta
$. Par ailleurs 
\begin{equation}\label{pluspetitram2}
\mbox{\sc  s}_{0} f (z) = a_0+
\int_0^{\infty}\widetilde{f}(\zeta)e^{-z \zeta} \, d\zeta = 
a_0 + \sum_{l=1}^{m} z^{\frac{m-l}{m}} 
\mbox{\sc  s}_{0} f_l (\dot z), \hspace{5mm} z \in P(B),
\end{equation}
et  pour tout $z \in P(B)$ et $n\geq 1$,
\begin{equation}\label{pluspetitram3}
 \Big|\mbox{\sc  s}_{0} f (z)  - \sum_{k=0}^{mn} 
  \frac{a_k}{z^\frac{k}{m}}\Big| \leq 
C e^{Br}  \frac{ n!  }{r^{n}  } \frac{\displaystyle  \sum_{i=0}^{m-1} |z|^\frac{i}{m}}{|z|^{n}(\Re
    (\dot z)-B)}, \hspace{5mm} C = \max_{1 \leq l \leq m} A_l .
\end{equation}
\end{Theorem}

\begin{proof} La preuve repose sur un analogue de 
la ``m\'ethode du vecteur cyclique''.
Ecrivons $f$ sous la forme
$$\displaystyle f(z) = a_0 + \sum_{l=1}^{m} z^{\frac{m-l}{m}} 
 f_l (\dot z), \hspace{5mm}\mbox{o\`u}
 \hspace{5mm}
 f_l(\dot z)=\sum_{j=1}^{+\infty}
\frac{a_{l,j}}{{\dot z}^{j}}, \, \, a_{l,j} = a_{l+m(j-1)}.$$
Puisque pour $k=0, 1, \cdots, m-1$,
$$ f(e^{2i\pi k}z) -a_0 =  \sum_{l=1}^{m}
 \omega^{-lk} z^\frac{m-l}{m}  f_l(\dot z), \hspace{5mm}\mbox{avec}
 \hspace{5mm} \omega = e^{\frac{2i\pi}{m}},$$ 
nous pouvons \'ecrire:
$$
A
\left(
\begin{array}{c}
 z^\frac{m-1}{m}  f_1(\dot z)\\
\\
 z^\frac{m-2}{m}  f_2(\dot z) \\
\\
\vdots \\
\\
  f_m(\dot z)
\end{array}
\right)=
\left(
\begin{array}{c}
 f^{[0]}(z)\\
\\
 f^{[1]} (z) \\
\\
\vdots \\
\\
 f^{[m-1]} (z)
\end{array}
\right), \hspace{5mm} \mbox{avec}  \hspace{5mm}  f^{[k]}(z) = f
(e^{2i\pi k} z) -a_0,
$$
o\`u 
$\displaystyle A = \left(
\begin{array}{ccc}
        & \vdots & \\
\cdots  &  A_{i,j} & \cdots  \\
        & \vdots & \\
\end{array}
\right), \hspace{2mm}  A_{i,j} = \omega^{-(i-1)j},
$ est une matrice $m \times m$ de Vandermonde inversible. Si
$\displaystyle A^{-1} = B=  \left(
\begin{array}{ccc}
        & \vdots & \\
\cdots  &  B_{i,j} & \cdots  \\
        & \vdots & \\
\end{array}
\right)
$, alors 
cela implique:
\begin{equation}\label{eqfl} 
 f_l (\dot z) = \frac{1}{z^\frac{m-l}{m}}\sum_{k=0}^{m-1} B_{l,(k+1)}
f^{[k]}(z),  \hspace{5mm} l=0, \cdots, m-1.
\end{equation}
La transform\'ee de Borel de $f^{[k]}(z)$ s'\'ecrivant sous la forme
$$\widetilde{f^{[k]}}(\zeta) = e^{-2i\pi k} \widetilde{f}(\zeta e^{-2i\pi k}),
$$
on en d\'eduit que chaque $\widetilde{f^{[k]}}(\zeta)$ se prolonge
analytiquement sur $\mathcal{D}_r^\star$ et
\begin{equation}\label{majsectfk}
\forall \zeta \in
\mathcal{D}_r^\star,  \, \, |\widetilde{f^{[k]}}(\zeta) \zeta^\frac{m-1}{m} |\leq A e
^{B|\zeta|}.
\end{equation}
De (\ref{eqfl}) on tire que
$$
\left\{
\begin{array}{l}
\displaystyle \widetilde{f_l}(\dot \zeta) =
\frac{\zeta^{-\frac{l}{m}}}{\Gamma(1-\frac{l}{m})} \ast \left(
  \sum_{k=0}^{m-1} B_{l,(k+1)} \widetilde{f^{[k]}}(\zeta)\right),
\hspace{5mm} l = 1, \cdots, m-1\\
\\
\displaystyle \widetilde{f_m}(\dot \zeta) =
  \sum_{k=0}^{m-1} B_{m,(k+1)} \widetilde{f^{[k]}}(\zeta).
\end{array}
\right.
$$
Ceci implique l'holomorphie de chaque $\widetilde{f_l}(\dot \zeta)$ sur
$\mathcal{B}_r$, et par ailleurs,
\begin{equation}\label{majsectfl}
\exists \, A_l >0, \, \forall \dot \zeta \in
\mathcal{B}_r,  \, \, |\widetilde{f_l}(\dot \zeta)  |\leq A_l e
^{B|\dot \zeta|}.
\end{equation}
Le th\'eor\`eme \ref{thmNevfin} s'applique alors \`a chacune des
s\'eries formelles Gevrey-1 $f_l(\dot z)$: pour  $\Re (\dot z)$ et $n\geq 1$,
$$ \Big|\mbox{\sc  s}_{0} f_l(\dot z) -\sum_{j=1}^{n}
  \frac{a_{l,j}}{{\dot z}^{j}}\Big| \leq 
A_le^{Br}  \frac{ n!  }{r^{n}} \frac{1}{|\dot z|^{n}(\Re
    (\dot z)-B)}
$$
o\`u 
$$\mbox{\sc  s}_{0} f_l(\dot z) =
\int_0^{\infty}\widetilde{f_l}(\zeta) e^{-\dot z \zeta} \, d\zeta.$$
En posant 
$$\mbox{\sc  s}_{0} f (z) = a_0 + \sum_{l=1}^{m} z^{\frac{m-l}{m}} 
\mbox{\sc  s}_{0} f_l (\dot z)$$
on d\'eduit de ce qui pr\'ec\'ede que pour tout $z \in P(B)$ et
$n\geq 1$,
$$ \Big|\mbox{\sc  s}_{0} f (z)  - \sum_{k=0}^{mn} 
  \frac{a_k}{z^\frac{k}{m}}\Big| \leq 
\max(A_l)e^{Br}  \frac{ n!  }{r^{n}  } \frac{\displaystyle  \sum_{i=0}^{m-1} |z|^\frac{i}{m}}{|z|^{n}(\Re
    (\dot z)-B)}.
$$
\end{proof}

\begin{rem}\label{precpluspetitram}
 La propri\'et\'e (\ref{majsectfl}) induit, par Cauchy, que  pour tout $j \geq 1$
\begin{equation}\label{plusfinraml}
|a_{l+m(j-1)}| \leq  A_l e^{Br} \frac{ j!  }{r^j}  \hspace{10mm} \Big(\mbox{car }
\displaystyle a_{l+m(j-1)} = \frac{d^j\widetilde{f_l}}{d\zeta ^j}(0) \Big).
\end{equation}
En pratique, on remplacera alors la majoration (\ref{pluspetitram3})
par une estimation de l'erreur de la forme: pour $\Re (\dot z)$ assez grand
et $n\geq 1$,
\begin{equation}\label{pluspetitram10}
 \Big|\mbox{\sc  s}_{0} f (z)  - \sum_{k=0}^{mn} 
  \frac{a_k}{z^\frac{k}{m}}\Big| \sim
\max_{1 \leq l \leq m}(|a_{l+mn}|)  
\frac{\displaystyle  \sum_{i=0}^{m-1} |z|^\frac{i}{m}}{|z|^{n}\Re
  (\dot z)}.
\end{equation}
Pour les m\^emes raisons que celles d\'evelopp\'ees \`a la remarque
\ref{precpluspetit}, la {\em sommation au plus petit terme} consistera
\`a choisir $n = \big[r|z|\big]$.
\end{rem}

\subsection{Un exemple}\label{exprat}

A titre d'exemple, qui nous servira \'egalement d'introduction \`a la section
\ref{FactGene},
nous allons consid\'erer la sommation de  Borel
d'une solution formelle de l'\'equation diff\'erentielle
\begin{equation}\label{ex1}
\frac{d^2 \Phi}{d x^2}=\frac{x^3-2x^2-3x+4}{x^2} \Phi.
\end{equation}

Suivant \cite{DelabRaso}:

\begin{Proposition}\label{newprop}
Soit $z(x) = \frac{2}{3}x^\frac{3}{2} -2x^\frac{1}{2}$. Il existe une unique
s\'erie formelle $\displaystyle \psi (z) \in
\mathbb{C}[[z^{-1/3}]]$,  de terme constant \'egal \`a
$1$, telle que 
\begin{equation}\label{Phi0Borel36}
\Phi  (x ) = \frac{e^{-z}}{z^\frac{1}{6}}  \psi (z) \,|_{z=z(x)}
\end{equation}
telle que $\Phi$ soit solution formelle de l'\'equation (\ref{ex1}). 
De plus  la transform\'ee de Borel de $\psi$ d\'efinit une fonction
analytique sur  le
 rev\^etement universel de $\mathbb{C}\backslash \{0, -2\}$ et est \`a
 croissance exponentielle d'ordre au plus 1 \`a l'infini. 
\end{Proposition}

La s\'erie formelle $\psi (z)$ de la proposition pr\'ec\'edente se
calcule \`a tout ordre:
\begin{equation}\label{cestunexemple}
\begin{array}{ll}
\displaystyle \psi  (z) & \displaystyle = \sum_{n=0}^{+\infty}
\frac{a_n}{z^{\frac{k}{3}}} \\
              & \displaystyle = 1-\left(\frac{128}{3}\right)^\frac{1}{3}
\frac{1}{z^\frac{1}{3}} + \left(\frac{2048}{9}\right)^\frac{1}{3}
\frac{1}{z^\frac{2}{3}} -\left(\frac{34328125}{373248}\right)^\frac{1}{3}
\frac{1}{z} + \cdots .
\end{array}
\end{equation}

Comme cons\'equence de la proposition \ref{newprop}, la s\'erie
formelle $\psi (z)$ rentre dans le cadre d'application du 
 th\'eor\`eme \ref{thmNevfinram}, avec $0 < r < 2$ : la 
s\'erie $\psi (z)$ est sommable de Borel pour $z \in P(B)$, $B>0$
assez grand,  et nous nous proposons
ici d'\'evaluer la somme 
 $\mbox{\sc  s}_{0}\psi (z)$. Celle-ci est de la forme
$$\mbox{\sc  s}_{0} \psi (z) = a_0 + \sum_{l=1}^{3} z^{\frac{3-l}{3}} 
\mbox{\sc  s}_{0} \psi_l (\dot z)$$
o\`u 
$$\psi_l(\dot z)=\sum_{j=1}^{+\infty} \frac{a_{l,j}}{{\dot z}^{j}},
\hspace{5mm} a_{l,j} = a_{l+3(j-1)}.$$

Pour l'illustration num\'erique qui suit, nous choisirons $z=12$ (cela
correspond \`a   $x=9$ dans la proposition \ref{newprop}).

\subsubsection{Sommation au plus petit terme} 
\begin{figure}[thp]
\begin{center}
\includegraphics[width=3.0in]{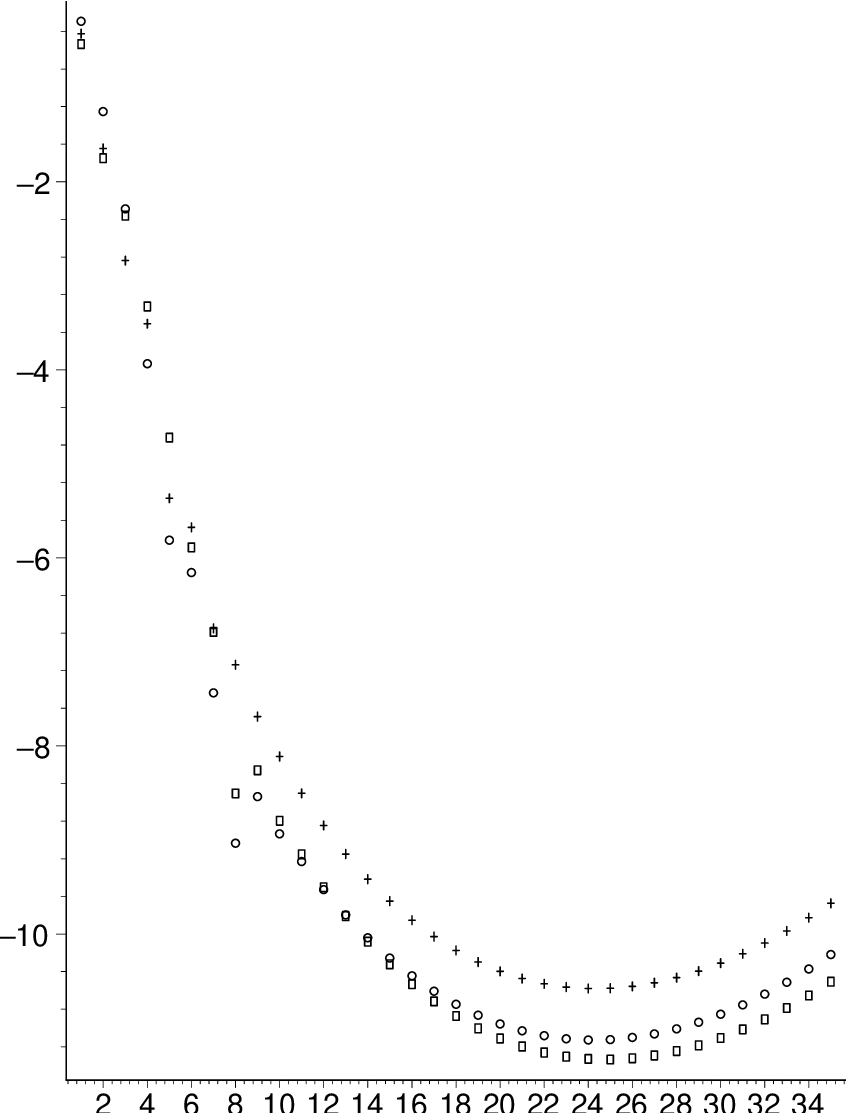} 
\caption{Nous avons repr\'esent\'e pour $z=12$ et 
 $j=1,\cdots, 35$, par $\Box$ les
  points de coordonn\'ees $[j, \log \left|\frac{a_{1,j}}{z^{j}}
  \right|]$,  par $\circ$ les
  points de coordonn\'ees $[j, \log \left|\frac{a_{2,j}}{z^{j}}
  \right|]$, par $+$ les
  points de coordonn\'ees $[j , \log \left|\frac{a_{3,j}}{z^{j}}
  \right|]$.
\label{fig:Som_Factorielle01}}
\end{center}
\end{figure}

Nous commen\c cons l'\'evaluation de 
$$\mbox{\sc  s}_{0} \psi (z),
 \hspace{5mm} z=12,$$
au moyen des sommes partielles $\displaystyle  \sum_{k=0}^{3n} 
  \frac{a_k}{z^\frac{k}{3}}$
par la sommation au plus petit terme comme expos\'e dans la remarque 
\ref{precpluspetitram}.\\
La figure \ref{fig:Som_Factorielle01} sugg\`ere de choisir $n=24$
comme troncation optimale, ce qui correspond au choix de 
$\displaystyle n = \sup_{0<r<2} \big[r|z|\big]$ (cf. 
 Remarque \ref{precpluspetitram}).  
 Le calcul donne:
$$\mbox{\sc  s}_{0} \psi (z) \simeq 0.26256292290 \pm 0.23 \times
10^{-9}.$$

\subsubsection{Sommation par s\'eries de factorielles}

Par la proposition
\ref{newprop} et le th\'eor\`eme  \ref{thmNevfinram}, 
les  transform\'ees
de Borel $\widetilde{ \psi_{(l)} }(\dot \zeta)$ des  $\psi_{(l)}
(\dot z)$ d\'efinissent des fonctions
 holomorphes dans le domaine $\mathcal{B}_r$ pour tout $0<r<2$,
mais \'egalement dans le domaine $\Delta_\lambda$ pour tout $0 < \lambda <
2/\ln(2)$, et sont \`a  croissance exponentielle d'ordre au
plus 1 dans ces domaines. En vertu du th\'eor\`eme
\ref{thmsomfact1ter}:

\noindent {\em 
Il existe $B>0$ tel que pour tout $0 < \lambda <
2/\ln(2)$ et tout $l=1,2,3$,   le d\'eveloppement
$$\lambda  \sum_{j=0}^{+\infty}
\frac{\Gamma(\lambda \dot z) \Gamma(j+1) b_{l,j}^{(\lambda)}}{\Gamma(\lambda
  \dot z+j+1)} $$
converge absolument pour $\Re (\dot z) > \max( B ,\frac{1}{\lambda})$ et sa somme repr\'esente
$\displaystyle \mbox{\sc s}_{0}\psi_l (\dot z)$, o\`u 
les coefficients $b_{l,j}^{(\lambda)}$ se
d\'eduisent des $a_{l,j}^{(\lambda)} = \lambda^{j-1} a_{l,j}$ par l'algorithme de Stirling.
}

\begin{table}[thp]
\begin{center}
\begin{tabular}{|c|l|c|}
\hline
Valeur de $n$ & Estimation de $\mbox{\sc  s}_{0} \psi (z)$ &
Estimation de l'erreur\\
\hline
10 & 0.262562935 & $0.20 \times 10^{-7}$\\
\hline
14 & 0.26256292301 & $0.22 \times 10^{-9}$\\
\hline
18 & 0.2625629228800 & $0.45 \times 10^{-11}$
\\
\hline
25 & 0.262562922877259 & $0.15 \times 10^{-13}$\\
\hline
33 & 0.262562922877250882 & $0.65 \times 10^{-16}$\\
\hline
40 & 0.2625629228772508441 & $0.2 \times 10^{-18}$\\
\hline
\end{tabular}
\caption{Calcul de $\mbox{\sc  s}_{0} \psi (z)$ par s\'eries de
  factorielles pour $z=12$ avec
  $\lambda = 2/\ln(2)$.
\label{table1}}
\end{center}
\end{table}

Evaluons \`a
pr\'esent la somme de Borel  
$$\mbox{\sc  s}_{0} \psi (z) = a_0 + \sum_{l=1}^{3} z^{\frac{3-l}{3}} 
\mbox{\sc  s}_{0} \psi_l (\dot z),  \hspace{5mm} z=12,$$
par l'utilisation des s\'eries de factorielles. 
On estime chacune des sommes de Borel $\mbox{\sc  s}_{0} \psi_{(l)}
(\dot z)$ au moyen des sommes partielles 
$\displaystyle  \lambda  \sum_{j=0}^n
\frac{\Gamma(\lambda \dot  z) \Gamma(j+1) b_{l,j}^{(\lambda)}}{\Gamma(\lambda
 \dot  z+j+1)} $. Les majorations  (\ref{remarquablelamb}) et
(\ref{pourbn}) am\'enent en pratique (pour  $\Re (\dot z)$ 
assez grand) \`a estimer l'erreur commise par
la relation 
\begin{equation}
 \Big|\mbox{\sc  s}_{0}  \psi_{(l)}(\dot z)  - \lambda  \sum_{j=0}^n
\frac{\Gamma(\lambda \dot  z) \Gamma(j+1) b_{l,j}^{(\lambda)}}{\Gamma(\lambda
 \dot  z+j+1)} \Big| \sim | b_{l,n+1}^{(\lambda)}| 
 \frac{|\Gamma (\lambda \dot  z)| \Gamma (n+1)}{\Re(\dot  z)
   |\Gamma(\lambda \dot  z+n+1)|}.
\end{equation}
En prenant  pour  $\lambda$  la borne sup
$2/\ln(2)$ des valeurs
th\'eoriquement permises, le calcul fournit la table \ref{table1}.

\begin{table}[thp]
\begin{center}
\begin{tabular}{|c|l|c|}
\hline
Valeur de $n$ & Estimation de $\mbox{\sc  s}_{0} \psi (z)$ &
Estimation de l'erreur\\
\hline
14 & 0.262562922891 & $0.24 \times 10^{-10}$\\
\hline
18 & 0.26256292287739 & $0.25 \times 10^{-12}$\\
\hline
\end{tabular}
\caption{Calcul de $\mbox{\sc  s}_{0} \psi (z)$ par s\'eries de
  factorielles pour $z=12$ avec
  $\lambda = 4$.
\label{table2}}
\end{center}
\end{table}

Ajoutons, sans tenter de l'expliquer, qu'on observe en pratique une acc\'el\'eration de la
convergence en prenant des valeurs de $\lambda$ au-del\`a des valeurs
th\'eoriquement permises. Ceci est illustr\'e par la table \ref{table2}.

\section{Sommation de Borel par s\'eries de factorielles des s\'eries
  de puissances fractionnaires}\label{FactGene}

L'exemple trait\'e au \S \ref{exprat} a illustr\'e comment la
m\'ethode de sommation d'une somme de Borel par les s\'eries de
factorielles pouvait \^etre adapt\'ee au cas d'une s\'erie de
puissances fractionnaires. 
Dans cette section nous allons pr\'esenter une variante
plus directe, qui peut \^etre vue
comme une extension de la m\'ethode de sommation par s\'eries de
factorielles.

\subsection{Sommation par s\'eries de factorielles g\'en\'eralis\'ees}

Nos hypoth\`eses seront les suivantes:  $\displaystyle f(z)=\sum_{n=0}^{+\infty}
\frac{a_n}{z^{\frac{n}{m}}} \in
\mathbb{C}[[z^{-\frac{1}{m}}]]_1$ est une s\'erie dont
la transform\'ee de Borel $\displaystyle \widetilde{f}(\zeta)=\sum_{n=1}^{+\infty}
\frac{a_n \zeta^{\frac{n}{m}-1} }{\Gamma\Big(\frac{n}{m}\Big)}  \in
\zeta^{-1+\frac{1}{m}}\mathbb{C}\{\zeta^\frac{1}{m}\} $ 
 se prolonge  analytiquement \`a l'ouvert
 $\Omega^\star$, et de plus:  
\begin{equation}\label{majOram}
\exists \, A >0, \,\exists \, B >0, \, \forall \zeta \in
\Omega^\star, \, \, 
|\widetilde{f}(\zeta) \zeta^\frac{m-1}{m} |\leq A e^{B|\zeta|}.
\end{equation}
 Comme $\mathcal{D}_{\ln(2)}^\star  \subset
 \Omega^\star$, on d\'eduit du 
th\'eor\`eme \ref{thmNevfinram} que 
la somme de Borel 
$$\displaystyle
\mbox{\sc  s}_{0} f(z) = 
a_0+\int_0^{+\infty}\widetilde{f}(\zeta)e^{-z \zeta} \, d\zeta$$
est bien d\'efinie pour $z \in P(B)$ et nous nous proposons de la calculer.\\

Pour $\zeta \in \Omega^\star$ nous pouvons \'ecrire $\widetilde{f}$ sous la forme
\begin{equation}\label{ideedesgl}
\left\{
\begin{array}{l}
 \displaystyle \widetilde{f}(\zeta) = \sum_{l=1}^{m}
(1-e^{-\zeta})^{\frac{l}{m}-1} \widetilde{g}_l(\dot \zeta)\\
\\
 \displaystyle \widetilde{g}_l(\dot \zeta) = 
 \left( \frac{\dot \zeta}{1-e^{-\dot \zeta}}\right)^{\frac{l}{m}-1}
\sum_{k=0}^{+\infty} \frac{a_{l+mk}
  {\dot \zeta}^k}{\Gamma\Big(\frac{l}{m}+k\Big)}  \in \mathbb{C}\{\dot \zeta\}
\end{array}
\right.
\end{equation}
et observons que 
pour tout $k \in \{0,1,\ldots,m-1\}$ nous avons:
$$ \widetilde{f}(e^{2i\pi k}\zeta)=\sum_{l=1}^{m} \omega^{kl}
\,(1-e^{-\zeta})^{\frac{l}{m}-1}
 \widetilde{g}_l(\dot \zeta) \hspace{5mm} \mbox{avec} \hspace{5mm} \omega=e^{\frac{2i\pi}{m}}.$$
Nous en tirons la relation~:
$$
A
\left(
\begin{array}{c}
\displaystyle  (1-e^{-\zeta})^{\frac{1}{m}-1} \widetilde{g}_1(\dot \zeta) \\
\vdots \\
\displaystyle (1-e^{-\zeta})^{\frac{m-1}{m}-1}
\widetilde{g}_{m-1}(\dot \zeta)\\
\\
\displaystyle \widetilde{g}_m(\dot \zeta)
\end{array}
\right)=
\left(
\begin{array}{c}
\displaystyle \widetilde{f}(e^{2i\pi }\zeta) \\
\vdots \\
\displaystyle  \widetilde{f}(e^{2i\pi (m-1)} \zeta)
\\
\displaystyle \widetilde{f}(\zeta)
\end{array}
\right).
$$
o\`u 
$\displaystyle A = \left(
\begin{array}{ccc}
        & \vdots & \\
\cdots  &  A_{i,j} & \cdots  \\
        & \vdots & \\
\end{array}
\right), \hspace{2mm}  A_{i,j} = \omega^{i(j-1)},
$ est une matrice $m \times m$ de Vandermonde inversible. En notant
$\displaystyle A^{-1} = B=  \left(
\begin{array}{ccc}
        & \vdots & \\
\cdots  &  B_{i,j} & \cdots  \\
        & \vdots & \\
\end{array}
\right)
$
on obtient que pour tout $l=1, \cdots, m$ et tout $\zeta \in \Omega^\star$,
$$\widetilde{g}_l(\dot \zeta) = (1-e^{-\zeta})^{1-\frac{l}{m}} 
\sum_{k=1}^m B_{l,k} \widetilde{f}(e^{2i\pi k} \zeta).$$
Cette propri\'et\'e, la d\'efinition m\^eme (\ref{ideedesgl}) des
$\widetilde{g}_l$, et les hypoth\`eses faites sur $\widetilde{f}$
montrent le lemme suivant:

\begin{Lemma}\label{croissanceexpdesgl}
Pour tout $l=1, \cdots, m$, $\widetilde{g}_l$ est holomorphe sur
$\Delta$ et  il existe $A_l
>0$ tel que, $\forall \dot \zeta \in \Delta$,  $\displaystyle
|\widetilde{g}_l(\dot \zeta)| \leq A_l e^{B|\dot \zeta|}$.
\end{Lemma}

Soit $D(1,1)^\star$  le disque ouvert \'epoint\'e de centre
$1$ et de rayon $1$. Notons
$\displaystyle 
\begin{array}{c}
D(1,1)_m^\star\\
\nu \downarrow \\
D(1,1)^\star
\end{array}$ 
le rev\^etement \`a $m$ feuillets de $D(1,1)^\star$. L'application
conforme $\dot  \zeta \in \Delta^\star \mapsto \dot s =e^{-\dot \zeta} \in
D(1,1)^\star$ se rel\`eve naturellement en une application conforme de
$ \Omega^\star$ sur $D(1,1)_m^\star$:
$$\begin{array}{ccc}
\zeta  \in \Omega^\star  & \longleftrightarrow & 
s =  e^{-\zeta} \in D(1,1)_m^\star \\
\pi \downarrow  & & \downarrow  \nu \\
\dot \zeta   \in \Delta^\star    &
\longleftrightarrow   &  
\dot s = e^{-\dot \zeta} \in D(1,1)^\star
\end{array}
$$
%%%%%%%%%%%%%%%%%
Posons alors, pour $s \in  D(1,1)_m^\star$:
$$\displaystyle
\Phi(s)=\widetilde{f}(\zeta),
\hspace{5mm} \phi_l (\dot s) = \widetilde{g_l} (\dot \zeta). $$
Suivant (\ref{ideedesgl})  l'application $\Phi$ se d\'ecompose sous la
forme
\begin{equation}\label{ram1}
 \displaystyle \Phi(s)=\sum_{l=1}^{m}
(1-s)^{\frac{l}{m}-1} \phi_l (\dot s)
\end{equation}
o\`u, comme dans la section \ref{factcla}, nous
pouvons \'ecrire, pour $\dot s \in D(1,1)$:
\begin{equation}\label{ram2}
\phi_l (\dot s)  =\sum_{j=0}^{+\infty}
 b_j^{(l)}(1-\dot s)^j.
\end{equation}
Par cons\'equent, $\Phi$ s'\'ecrit  sous la forme:
\begin{equation}\label{ram3}
\begin{array}{l}
\displaystyle \forall \, s \in  D(1,1)_m^\star, \, \, \,
\Phi(s)=\sum_{n=1}^{+\infty} d_n (1-s)^{\frac{n}{m}-1},\\
\\
\displaystyle \mbox{avec } \, \,  \forall l \in \{1,\ldots,m\}, \, \forall j \in
\mathbb{N}, \quad d_{l+mj}= b_{j}^{(l)}.
\end{array}
\end{equation}
Formellement, nous pouvons \'ecrire la somme de Borel de $f$ sous la forme:
$$ \begin{array}{ll}
\mbox{\sc  s}_{0} f(z)  &  \displaystyle =  a_0+\int_0^{+\infty}\widetilde{f}(\zeta)e^{-z
  \zeta} \, d\zeta = a_0+\int_0^{+\infty}
\sum_{n=1}^{+\infty} d_n (1-e^{-\zeta})^{\frac{n}{m}-1}e^{-z \zeta} \,
d\zeta, \\
\\
  & \displaystyle = a_0 + \sum_{n=1}^{+\infty}d_n\int_0^1
(1-s)^{\frac{n}{m}-1}s^{z-1} \,ds = 
a_0+ \sum_{n=1}^{+\infty} \frac{\Gamma\Big(\frac{n}{m}\Big)\Gamma(z) d_n}{\Gamma\Big(z+\frac{n}{m}\Big)}.
\end{array}$$

Cette derni\`ere expression de $\mbox{\sc  s}_{0} f(z)$ est
justifi\'ee par la 
g\'en\'eralisation suivante du th\'eor\`eme \ref{thmsomfact1}:

\begin{Proposition}\label{thmsomfact2}
Soit $\displaystyle f(z)=\sum_{n=0}^{+\infty}
\frac{a_n}{z^{\frac{n}{m}}} \in
\mathbb{C}[[z^{-\frac{1}{m}}]]_1$. On suppose que 
la transform\'ee de Borel $\displaystyle \widetilde{f}(\zeta)=\sum_{n=1}^{+\infty}
\frac{a_n \zeta^{\frac{n}{m}-1} }{\Gamma\Big(\frac{n}{m}\Big)}  \in
\zeta^{-1+\frac{1}{m}}\mathbb{C}\{\zeta^\frac{1}{m}\} $ 
 se prolonge  analytiquement \`a l'ouvert
 $\Omega^\star$ et  que
$$
\exists \, A >0, \,\exists \, B >0, \, \forall \zeta \in
\Omega^\star, \, \, 
|\widetilde{f}(\zeta) \zeta^\frac{m-1}{m} |\leq A e^{B|\zeta|}.
$$
Alors la s\'erie de factorielles g\'en\'eralis\'ee
$$\displaystyle a_0+\sum_{n=1}^{+\infty}
\frac{\Gamma\Big(\frac{n}{m}\Big)\Gamma(z)
  d_n}{\Gamma\Big(z+\frac{n}{m}\Big)}$$
converge absolument pour $z \in P(\max(B,1))$
et repr\'esente la somme de Borel $\mbox{\sc  s}_{0} f(z)$ dans cet ouvert.
\end{Proposition}

\begin{proof}
La preuve s'appuiera sur deux lemmes pr\'eparatoires.

\begin{Lemma}\label{eqBeta}
$ \displaystyle \forall z \in \mathbb{C}\backslash \mathbb{R}^-, \quad
\frac{\Gamma(z)\Gamma(\frac{n}{m})}{\Gamma(z+\frac{n}{m})}
\stackrel{\sim}{_{n \to \infty}} \frac{\Gamma(z)}{\big(\frac{n}{m}\big)^z}.$
\end{Lemma}

\begin{proof}
Nous savons par la formule de Stirling que pour $z \in \mathbb{C}\backslash{\mathbb{R}^-}$,
\break
 $ \displaystyle \Gamma(z) \stackrel{\sim}{_{\substack{|z| \to
    +\infty \\ z \in \mathbb{C}\backslash \mathbb{R}^-}}}
\sqrt{2\pi}z^{z-\frac{1}{2}}e^{-z}$, d'o\`u~:
$$ \left\{
\begin{array}{l}
\displaystyle \Gamma(\frac{n}{m}) \stackrel{\sim}{_{n \to +\infty}}
\sqrt{2\pi}\big(\frac{n}{m}\big)^{\frac{n}{m}-\frac{1}{2}}e^{-\frac{n}{m}}\\
\\
\displaystyle \frac{1}{\Gamma(z+\frac{n}{m})} \stackrel{\sim}{_{n \to +\infty}} \frac{e
  ^{z+\frac{n}{m}}}{\sqrt{2\pi}\big(z+\frac{n}{m}\big)^{z+\frac{n}{m}-\frac{1}{2}}} = \frac{e
  ^{\frac{n}{m}}}{\sqrt{2\pi}\big(z+\frac{n}{m}\big)^{z-\frac{1}{2}}}
\, \frac{e
  ^{z}}{\big(z+\frac{n}{m}\big)^{\frac{n}{m}}}.
\end{array}
\right.$$
Or, nous avons les \'equivalences suivantes~:
$$ \left\{
\begin{array}{l}
\displaystyle \big(z+\frac{n}{m}\big)^{z-\frac{1}{2}} \stackrel{\sim}{_{n \to +\infty}}
\big(\frac{n}{m}\big)^{z-\frac{1}{2}}\\
\\
\displaystyle \big(z+\frac{n}{m}\big)^{-\frac{n}{m}} \stackrel{\sim}{_{n \to
  +\infty}} \big(\frac{n}{m}\big)^{-\frac{n}{m}}e^{-z},
\end{array}
\right.$$
ce qui nous donne l'\'equivalence souhait\'ee. 
\end{proof} 

\begin{Lemma}\label{CVabs}
Consid\'erons la fonction $\displaystyle \Phi(s)= \sum_{n=1}^{+\infty} d_n (1-s)^{\frac{n}{m}-1}$.\\
Sous les hypoth\`eses de la proposition \ref{thmsomfact2},
la s\'erie $\displaystyle \sum_{n=1}^{+\infty}
\frac{|d_n|}{\big(\frac{n}{m}\big)^{C}}$ converge pour tout $C>\max(B,1)$.
\end{Lemma}

\begin{proof}
Par le lemme \ref{croissanceexpdesgl} et le  lemme \ref{lemmecruc}
nous pouvons d\'eduire que les coefficients $b_j^{(l)}$ d\'efinis par
(\ref{ram2}) v\'erifient:
$$ \forall C> \max(B,1), \hspace{5mm} \sum_{j=1}^{+\infty}  \frac{|b_j^{(l)}|}{j^{C}}
< +\infty.$$
A fortiori, pour tout $l=1, \cdots, m$,
$$ \forall C> \max(B,1) , \hspace{5mm} \sum_{j=0}^{+\infty}  \frac{|b_j^{(l)}|}{
\left(j+\frac{l}{m}\right)^{C}}
< +\infty.$$
Par suite, par la d\'efinition (\ref{ram3}) des coefficients $d_n$ et
par sommation finie sur $l$,
$$ \forall C> \max(B,1), \hspace{5mm} \sum_{l=1}^m
\sum_{j=0}^{+\infty}  \frac{ |d_{l+mj}| }{ \left(j+\frac{l}{m}\right)^{C} }
< +\infty.$$
Ceci fournit {\sl in fine} la relation:
$\displaystyle  \forall C> \max(B,1) , \hspace{5mm} \sum_{n=1}^{+\infty} \frac{
  |d_n| }{ \left(\frac{n}{m}\right)^{C} }
< +\infty$.
\end{proof}

Nous revenons maintenant \`a la preuve de la proposition
\ref{thmsomfact2} proprement dite.\\
 Pour ce qui est de la convergence absolue de la s\'erie
$\displaystyle a_0+\sum_{n=1}^{+\infty}
\frac{\Gamma\Big(\frac{n}{m}\Big)\Gamma(z)
  d_n}{\Gamma\Big(z+\frac{n}{m}\Big)}$ 
pour $z \in P(\max(B,1))$, comme par le lemme \ref{eqBeta}:
$$\displaystyle  d_n\frac{\Gamma(z)\Gamma(\frac{n}{m})}{\Gamma(z+\frac{n}{m})} \stackrel{\sim}{_{n
  \to +\infty}} d_n\Gamma(z)\Big(\frac{n}{m}\Big)^{-z},$$ 
il suffit donc de voir
  que pour $z \in P(\max(B,1))$ la s\'erie  $\displaystyle
\sum_{n=1}^{+\infty} |d_n|\Big(\frac{n}{m}\Big)^{-\Re(\dot z)}$
converge, ce qui est une cons\'equence du lemme \ref{CVabs}.\\
Pour voir que la s\'erie $\displaystyle a_0+\sum_{n=1}^{+\infty}
\frac{\Gamma\Big(\frac{n}{m}\Big)\Gamma(z)
  d_n}{\Gamma\Big(z+\frac{n}{m}\Big)}$ repr\'esente bien la somme de Borel $\mbox{\sc  s}_{0} f(z)$, il
  s'agit de montrer que dans l'expression
$$ \displaystyle \sum_{n=1}^{+\infty} |d_n|\int_0^1 s^{z-1}(1-s)^{\frac{n}{m}-1}ds$$
nous pouvons permuter $\sum$ et $\int$.
Il suffit pour cela de montrer que la fonction $ \displaystyle \sum_{n=1}^{+\infty}
|d_n||s^{z-1}|(1-s)^{\frac{n}{m}-1}$
est int\'egrable sur $[0,1]$.
Or, nous avons les \'egalit\'es suivantes (en posant $C=\Re(\dot z)>\max(B,1)$):
$$ \begin{array}{lll}
\displaystyle \sum_{n=1}^{+\infty}
|d_n|\int_0^1(1-s)^{\frac{n}{m}-1}|s^{z-1}|ds & =
& \displaystyle \sum_{n=1}^{+\infty}
|d_n|\int_0^1(1-s)^{\frac{n}{m}-1}s^{C-1}ds\\ 
\\ 
 & = & \displaystyle \sum_{n=1}^{+\infty} |d_n|
\frac{\Gamma\Big(\frac{n}{m}\Big)\Gamma(C)}{\Gamma\Big(C+\frac{n}{m}\Big)}. 
\end{array}$$
Or cette derni\`ere s\'erie converge comme nous l'avons d\'emontr\'e au 
point pr\'ec\'edent. Ceci ach\`eve la d\'emonstration.
\end{proof}
 
Il nous reste pour terminer \`a d\'eduire les coefficients $d_n$ des $a_n$. Tout ce
que nous avons \`a faire est de calculer la d\'ecomposition donn\'ee
par la proposition \ref{thmsomfact2} 
pour  $\displaystyle \frac{1}{z^{r}}$, $\displaystyle r >0$.
Pour $z \in P(0)$, nous avons
$$\frac{1}{z^{r}} = \frac{1}{\Gamma(r)}\int_0^{+\infty}
e^{-uz}u^{r-1} du$$
de sorte que, avec $\displaystyle u = -\ln (s )$,
$\displaystyle \frac{1}{z^{r}} = \frac{1}{\Gamma(r)}\int_0^1
s^{z-1} \left(-\ln (s)\right)^{r-1} ds$.
Pour $s \in ]0,1[$, nous pouvons \'ecrire
$\displaystyle -\frac{\ln (s)}{1-s}  = 
\sum_{j=0}^{+\infty}  \frac{j !}{j+1} \frac{(1-s)^j}{j !}$.
La s\'erie de Taylor,
$$ \left(- \frac{\ln (s)}{1-s} \right)^{r-1} = 
\sum_{j=0}^{+\infty} c_{r, j} (1-s)^j,$$
se d\'eduit alors de la formule de Faa di Bruno (\cite{Comtet}), et nous obtenons~:
$$c_{r, 0} =1, \hspace{5mm}
 c_{r, j} =  \frac{1}{j !} \sum_{1 \leq p \leq j } \frac{\Gamma (r)}{\Gamma (r-p)}
  B_{j,p} (\frac{1 !}{2}, \frac{2 !}{3}, \cdots, \frac{l !}{l+1}, 
  \cdots ),  \, \, j \geq 1, $$ 
o\`u les $B_{j,p}$ d\'esignent les polyn\^omes de Bell exponentiels
partiels  (\cite{Comtet}). En permutant $\sum$ et $\int$ (licite par un
calcul identique \`a celui effectu\'e dans la d\'emonstration de la
proposition \ref{thmsomfact2}), nous en
d\'eduisons que
$$\frac{1}{z^{r}} = \sum_{j=0}^{+\infty} \frac{ c_{r,j}}{\Gamma(r) }\int_0^1
s^{z-1} (1-s)^{r+j -1} ds =
\sum_{j=0}^{+\infty} \frac{ c_{r,j}}{\Gamma(r) }\beta(r+j,z)
$$
$$\frac{1}{z^{r}} = \sum_{j=0}^{+\infty} \frac{ c_{r,j}}{\Gamma(r) } 
\frac{\Gamma(r+j) \Gamma (z)}{\Gamma(r+j +z)}.
$$
Nous avons en particulier:
$$\frac{1}{z^{r}} = 
\frac{\Gamma (z)}{\Gamma(r+z)} 
+
\sum_{j=1}^{+\infty}  d_{r,j}
\frac{ \Gamma (z)}{\Gamma(r+j +z) }
$$
avec
$$
d_{r,j} = \left( \sum_{1 \leq p \leq j } 
\frac{  B_{j,p} (\frac{1 !}{2}, \frac{2 !}{3}, \cdots, \frac{l !}{l+1}, 
  \cdots )}{\Gamma (r-p)}
  \right) \frac{\Gamma(r+j)}{j !}.
$$
En particulier, pour $m \in \mathbb{N}^\star$ et $l \in
\mathbb{N}^\star$,
$$\frac{1}{z^{\frac{l}{m}}} = 
\frac{\Gamma (z)}{\Gamma(z+\frac{l}{m})} 
+
\sum_{j=1}^{+\infty}  d_{\frac{l}{m},j}
\frac{ \Gamma (z)}{\Gamma(z +\frac{l+jm}{m}) }.
$$
Nous en d\'eduisons alors facilement le r\'esultat qui suit:

\begin{Proposition}\label{lemcalram}
Dans la proposition \ref{thmsomfact2}, nous avons, pour $n \in \mathbb{N}^\star$,
$$d_n = \frac{1}{\Gamma(\frac{n}{m})}\left( a_n + 
\renewcommand{\arraystretch}{0.5}
\begin{array}[t]{c}
\sum \\
{\scriptstyle j \geq 1, \, l \geq 1 }\\
{\scriptstyle l+jm = n}
\end{array}
\renewcommand{\arraystretch}{1}
d_{\frac{l}{m},j}. a_l \right)
$$
o\`u les $d_{r,j}$ sont d\'efinis par
\begin{equation}\label{ram4}
d_{r,j} = \left( \sum_{1 \leq p \leq j } 
\frac{  B_{j,p} (\frac{1 !}{2}, \frac{2 !}{3}, \cdots, \frac{l !}{l+1}, 
  \cdots )}{\Gamma (r-p)}
  \right) \frac{\Gamma(r+j)}{j !},
\end{equation}
les $B_{j,p}$ d\'esignant les polyn\^omes de Bell exponentiels
partiels.
\end{Proposition}

La proposition \ref{thmsomfact2} induit 
le r\'esultat  suivant:

\begin{Theorem}\label{thmsomfact2bis}
Soit $\displaystyle f(z)=\sum_{n=0}^{+\infty}
\frac{a_n}{z^{\frac{n}{m}}} \in
\mathbb{C}[[z^{-\frac{1}{m}}]]_1$. On suppose qu'il existe $\lambda
>0$ tel que
la transform\'ee de Borel $\displaystyle \widetilde{f}(\zeta)=\sum_{n=1}^{+\infty}
\frac{a_n \zeta^{\frac{n}{m}-1} }{\Gamma\Big(\frac{n}{m}\Big)}  \in
\zeta^{-1+\frac{1}{m}}\mathbb{C}\{\zeta^\frac{1}{m}\} $ 
 se prolonge  analytiquement \`a l'ouvert
 $\Omega_\lambda^\star$, et  que
$$
\exists \, A >0, \,\exists \, B >0, \, \forall \zeta \in
\Omega_\lambda^\star, \, \, 
|\widetilde{f}(\zeta) \zeta^\frac{m-1}{m} |\leq A e^{B|\zeta|}.
$$
Alors la s\'erie de factorielles g\'en\'eralis\'ee
$$a_0+\lambda\sum_{n=1}^{+\infty}
\frac{\Gamma\Big(\frac{n}{m}\Big)\Gamma(\lambda z)
  d_n^{(\lambda)}}{\Gamma\Big(\lambda z+\frac{n}{m}\Big)},
$$
o\`u les $d_n^{(\lambda)}$ se d\'eduisent des $a_n^{(\lambda)} =
\lambda^{\frac{n}{m}-1} a_n$ par la proposition \ref{lemcalram},
converge absolument pour $z \in P(\max(B,1/\lambda))$
et repr\'esente la somme de Borel $\mbox{\sc  s}_{0} f(z)$ dans cet ouvert.
\end{Theorem}

\begin{proof}
Elle est similaire \`a celle du th\'eor\`eme \ref{thmsomfact1ter}.
\end{proof}

\subsection{Exemple 1}

Nous reprenons l'exemple de la sous-section \ref{exprat}. Nous
estimons
 la somme de Borel $\mbox{\sc  s}_{0} \psi(z)$
pour $z=12$ au moyen de la s\'erie de factorielles g\'en\'eralis\'ee
tronqu\'ee
$$a_0+\lambda\sum_{k=1}^{N}
\frac{\Gamma\Big(\frac{k}{3}\Big)\Gamma(\lambda z)
  d_k^{(\lambda)}}{\Gamma\Big(\lambda z+\frac{k}{3}\Big)},
$$
avec $\lambda=2/\ln(2)$ et $N=n/3$. La comparaison, d\'etaill\'ee par
la table \ref{table3},  est faite avec la
``valeur exacte'' (voir table \ref{table1}),
$$
\mbox{``valeur exacte''} = 0.2625629228772508441 \pm 0.2 \times 10^{-18}.
$$

\begin{table}[thp]
\begin{center}
\begin{tabular}{|c|l|c|}
\hline
Valeur de $n=N/3$ & Estimation de $\mbox{\sc  s}_{0} \psi (z)$ &
Erreur\\
\hline
10 & 0.262562936 & $0.13 \times 10^{-7}$\\
\hline
18 & 0.2625629228786 & $0.13 \times 10^{-11}$
\\
\hline
25 & 0.2625629228772537 & $0.29 \times 10^{-14}$\\
\hline
\end{tabular}
\caption{Calcul de $\mbox{\sc  s}_{0} \psi (z)$ par s\'eries de
  factorielles g\'en\'eralis\'ees pour $z=12$ avec
  $\lambda = 2/\ln(2)$.
\label{table3}}
\end{center}
\end{table}

\subsection{Exemple 2}

Consid\'erons \`a pr\'esent la somme de Borel
$$ \mbox{\sc  s}_{0} f(z) = \int_0^{+\infty}
\left(1+\zeta^{1/2}\right)^{1/2} e^{-z\zeta} \, d\zeta$$
de la s\'erie Gevrey
$$f(z) = \sum_{k =0}^\infty (-1)^{k+1}
\frac{\Gamma(\frac{k}{2}+1)\Gamma(k-\frac{1}{2})}{2\sqrt{\pi} \Gamma(k+1)}\frac{1}{z^{1+k/2}}$$
qui rentre dans le cadre de la proposition \ref{MagGevrey}, mais pas
dans celui des th\'eor\`emes \ref{thmNevfinram} et
\ref{thmsomfact2bis}, du fait de la singularit\'e en $\zeta =
e^{2i\pi}$ pour la transform\'ee de Borel. Un calcul direct montre que
$$\mbox{\sc  s}_{0} f(5) = 0.2357006$$
\`a $10^{-7}$ pr\`es. Le calcul \`a $10^{-6}$ pr\`es par s\'eries de
  factorielles g\'en\'eralis\'ees tronqu\'eees avec $\lambda=1$,
$\displaystyle \sum_{k=1}^{N}
\frac{\Gamma\Big(\frac{k}{2}\Big)\Gamma( z)
  d_k} {\Gamma\Big( z+\frac{k}{2}\Big)}
$, donne la table \ref{table4}: comme on pouvait le pr\'evoir, la
s\'erie de factorielles g\'en\'eralis\'ee ne converge pas
vers la somme de Borel $\mbox{\sc  s}_{0} f(5)$.

\begin{table}[thp]
\begin{center}
\begin{tabular}{|c|c|}
\hline
Valeur de $N$ & Estimation par s\'eries de factorielles \\
\hline
10 & 0.235584 \\
\hline
100 & 0.159338 \\
\hline
\end{tabular}
\caption{
\label{table4}}
\end{center}
\end{table} 

Pour se tirer d'affaire on peut utiliser la remarque \ref{somplusram}:
en prenant $\displaystyle \theta = \frac{\pi}{3}$ (par exemple), la
somme de Borel  $\mbox{\sc  s}_{\theta} f (z)$  d\'efinit un
prolongement analytique de $\mbox{\sc  s}_{0} f (z)$ pour
$\displaystyle \arg (z) \in ]-\frac{5\pi}{6}, \frac{\pi}{6}[$,
$|z|>0$. Donc $\mbox{\sc  s}_{\theta} f (5) = \mbox{\sc  s}_{0} f (5)$. Or
$$
\mbox{\sc  s}_{\theta} f (z) = \mbox{\sc  s}_{0} f_\theta (ze^{i\theta})
$$
o\`u
$$f_\theta (z) = f(ze^{-i\theta})  = \sum_{k =0}^\infty (-1)^{k+1}
\frac{\Gamma(\frac{k}{2}+1)\Gamma(k-\frac{1}{2})}{2\sqrt{\pi}
  \Gamma(k+1)}\frac{e^{i\frac{\pi}{3}(1+\frac{k}{2})}}{z^{1+k/2}}.$$
Les th\'eor\`emes \ref{thmNevfinram} et
\ref{thmsomfact2bis} s'appliquent \`a $f_\theta$, sous r\'eserve de
prendre $\lambda$ tel que l'image de $\Delta_\lambda$ par la rotation
de centre $0$ et d'angle $\theta$ ne contienne pas $1$. On peut prendre
$\lambda=0.6$ par exemple. Ceci permet
l'\'evaluation de $\mbox{\sc  s}_{0} f_\theta (z)$ par la s\'erie de
factorielles g\'en\'eralis\'ee associ\'ee
$\displaystyle \sum_{k=1}^{\infty}
\frac{\Gamma\Big(\frac{k}{2}\Big)\Gamma( \lambda z)
  d_k^{(\lambda)}(\theta) }{\Gamma\Big(\lambda z+\frac{k}{2}\Big)}
$. On estime alors $\mbox{\sc  s}_{0} f (5)$ en \'evaluant les sommes
partielles
$\displaystyle \sum_{k=1}^{N}
\frac{\Gamma\Big(\frac{k}{2}\Big)\Gamma( \lambda z)
  d_k^{(\lambda)}(\theta) }{\Gamma\Big(\lambda z+\frac{k}{2}\Big)}
$ pour $z = 5 e^{i \theta}$. Le r\'esultat est illustr\'e par la table
\ref{table5}, la convergence \'etant tr\`es lente.

\begin{table}[thp]
\begin{center}
\begin{tabular}{|c|c|c|}
\hline
Valeur de $N$ & Estimation par s\'eries de factorielles & |erreur|\\
\hline
50 & $0.2356902 +0.50\times10^{-5} i$ & $0.12\times10^{-4}$\\
\hline
150 & $0.2357024-  0.25 \times10^{-6} i $ & $0.1 \times10^{-5}$\\
\hline
\end{tabular}
\caption{
\label{table5}}
\end{center}
\end{table}

%%%%%%%%%%%%%%%%%%%%%%%%%%%%%%%%%%%%%%%%%%%%%

\end{document}